\documentclass[11pt]{article}
\usepackage{amsfonts}
\usepackage{color}
\newtheorem{theo}{Theorem}[section]
\newtheorem{lem}[theo]{Lemma}
\newtheorem{cor}[theo]{Corollary}

\newtheorem{defi}{Definition}[section]

\newcommand{\mysection}[1]{\section{#1} \setcounter{equation}{0}}
\newcommand{\proof}{{\sc Proof.} \quad}
\newcommand{\proofc}{{\sc Proof} \ }
\newcommand{\be}{\begin{equation} \label}
\newcommand{\ee}{\end{equation}}
\newcommand{\bea}{\begin{eqnarray}\label}
\newcommand{\eea}{\end{eqnarray}}
\newcommand{\bas}{\begin{eqnarray*}}
\newcommand{\eas}{\end{eqnarray*}}
\newcommand{\bit}{\begin{itemize}}
\newcommand{\eit}{\end{itemize}}
\newcommand{\qed}{\hfill$\Box$ \vskip.2cm}
\newcommand{\nn}{\nonumber}
\newcommand{\R}{\mathbb{R}}
\newcommand{\N}{\mathbb{N}}
\newcommand{\pO}{\partial\Omega}

\newcommand{\eps}{\varepsilon}

\newcommand{\wto}{\rightharpoonup}
\newcommand{\wsto}{\stackrel{\star}{\rightharpoonup}}
\newcommand{\hra}{\hookrightarrow}
\newcommand{\io}{\int_\Omega}

\newcommand{\abs}{\\[5pt]}
\newcommand{\Abs}{\\[5mm]}

\newcommand{\proj}{{\cal P}}
\newcommand{\neps}{n_\eps}
\newcommand{\ceps}{c_\eps}
\newcommand{\ueps}{u_\eps}
\newcommand{\Peps}{P_\eps}
\newcommand{\tme}{T_{max,\eps}}
\newcommand{\kd}{k_D}
\newcommand{\nzb}{\overline{n_0}}
\newcommand{\onz}{\overline{n_0}}
\begin{document}
\enlargethispage{10mm}
\title{Boundedness and large time behavior in a three-dimensional chemotaxis-Stokes system 
with nonlinear diffusion \\
and general sensitivity}
\author{
Michael Winkler\footnote{michael.winkler@math.uni-paderborn.de}\\
{\small Institut f\"ur Mathematik, Universit\"at Paderborn,}\\
{\small 33098 Paderborn, Germany} }
\date{}
\maketitle
\begin{abstract}
\noindent 
  We consider the chemotaxis-fluid system 
    \be{00}
    \left\{ \begin{array}{rcll}
    n_t + u\cdot\nabla n &=& \nabla \cdot (D(n)\nabla n) - \nabla \cdot (nS(x,n,c)\cdot \nabla c),\\[1mm]
    c_t + u\cdot\nabla c &=& \Delta c-nf(c),    \\[1mm]
    u_t + \nabla P &=& \Delta u   + n \nabla \phi,  \\[1mm]
    \nabla \cdot u &=& 0,
    \end{array} \right.
  \ee
  in a bounded convex domain $\Omega\subset \R^3$ with smooth boundary, where $\phi\in W^{1,\infty}(\Omega)$ 
  and $D, f$ and $S$ are given functions with values in $[0,\infty), [0,\infty)$ and $\R^{3\times 3}$, respectively.\abs
  In the existing literature, the derivation of results on global existence and qualitative behavior 
  essentially relies on the use of
  energy-type functionals which seem to be available only in special situations, 
  necessarily requiring the matrix-valued
  $S$ to actually reduce to a scalar function of $c$ which, along with $f$, in addition should satisfy
  certain quite restrictive structural conditions.\abs
  The present work presents a novel a priori estimation method 
  which allows for removing any such additional hypothesis:
  Besides appropriate smoothness assumptions, in this paper it is only required that $f$ 
  is locally bounded in $[0,\infty)$, that 
  $S$ is bounded in $\Omega\times [0,\infty)^2$, 
  and that $D(n)\ge \kd n^{m-1}$ for all $n\ge 0$ with some $\kd>0$ and some
  \bas
	m>\frac{7}{6}.
  \eas
  It is shown that then for all reasonably regular initial data, a corresponding initial-boundary value problem for (\ref{00})
  possesses a globally defined weak solution.\abs
  The method introduced here is efficient enough to moreover provide global boundedness of all solutions
  thereby obtained in that, inter alia, $n\in L^\infty(\Omega\times (0,\infty))$. 
  Building on this boundedness property, it can finally even be 
  proved that in the large time limit, any such solution approaches 
  the spatially homogeneous equilibrium $(\onz,0,0)$ in an appropriate sense, where $\onz:=\frac{1}{|\Omega|} \io n_0$,
  provided that merely $n_0\not\equiv 0$ and $f>0$ on $(0,\infty)$.
  To the best of our knowledge, these are the first results on boundedness and asymptotics of large-data solutions in a 
  three-dimensional chemotaxis-fluid system of type (\ref{00}).\abs
\end{abstract}
\newpage
\section{Introduction}\label{intro}
We consider the chemotaxis-Stokes system
\be{0}
    \left\{ \begin{array}{rcll}
    n_t + u\cdot\nabla n &=& \nabla \cdot \Big(D(n)\nabla n\Big)  - \nabla \cdot \Big(nS(x,n,c)\cdot \nabla c\Big), 
	\qquad & x\in\Omega, \ t>0,\\[1mm]
    c_t + u\cdot\nabla c &=& \Delta c-nf(c), \qquad & x\in\Omega, \ t>0,   \\[1mm]
     u_t   &=& \Delta u - \nabla P+ n \nabla \phi, \qquad & x\in\Omega, \ t>0, \\[1mm]
    \nabla \cdot u &=& 0, \qquad & x\in\Omega, \ t>0,
    \end{array} \right.
\ee
in a bounded domain $\Omega \subset \R^N$, where the main focus of this work will be on the case $N=3$.
Systems of this type arise in the modeling of populations of aerobic bacteria when suspended into sessile drops of water
(\cite{goldstein2004}, \cite{goldstein2005}).
In this setting, $n=n(x,t)$ and $c=c(x,t)$ denote the density of the cell population and the oxygen concentration,
respectively, and $u=u(x,t)$ and $P=P(x,t)$ represent the fluid velocity and the associated pressure.
The essential modeling hypotheses underlying (\ref{0}) are that cell movement is partially directed by gradients
of the chemical which they consume, that convection transports both cells and oxygen, and that the presence of 
bacteria, which are slightly heavier than water, influences the fluid motion through buoyant forces in an external
gravitational potential $\phi$.
The additional assumption that the fluid flow be comparatively slow is reflected in the fact that in (\ref{0})
its evolution is described by the Stokes equations rather than the full Navier-Stokes system
(\cite{lorz}).
Related mechanisms of chemotaxis-fluid interaction also arise in different biological contexts such as biomixing-based
fertilization strategies of certain benthic invertebrates (\cite{kiselev_ryzhik_CPDE2012}, \cite{kiselev_ryzhik_JMP2012}).\Abs
{\bf Approaches based on a natural energy functional.} \quad
According to the model specification underlying the numerical simulations in the original work \cite{goldstein2005},
analytical studies in the existing mathematical literature concentrate on the particular version of (\ref{0}) 
obtained on considering
\be{1}
    n_t + u\cdot\nabla n = \nabla \cdot \Big(n^{m-1}\nabla n\Big)  - \nabla \cdot \Big(n\chi(c) \nabla c\Big), 
	\qquad x\in\Omega, \ t>0,
\ee
as the first equation therein,
with $m\ge 1$ and a scalar chemotactic sensitivity function $\chi:\ [0,\infty)\to\R$.
Here an essential step forward in the analysis was marked by the observation that under 
suitable structural assumptions linking $\chi$ to the oxygen consumption rate $f$, 
this class of versions of (\ref{0}) admits for certain
natural quasi-Lyapunov functionals which involve the logarithmic entropy $\io n\ln n$.
Indeed, when tracking the time evolution of the latter, the appearing crucial cross-diffusion-related integral 
$\io \chi(c)\nabla c \cdot \nabla n$ can precisely be cancelled upon adding the result of 
a suitable testing procedure in the second equation
in (\ref{0}), where further integrals arising during the latter can be controlled conveniently 	
under certain conditions on the relationship between $\chi$ and $f$.
In the prototypical case $m=1, \chi\equiv 1$ and $f(c)=c$, for instance, this 
gives rise to an inequality of the form
\be{energy}
	\frac{d}{dt} \bigg\{ \io n\ln n +\frac{1}{2} \io \frac{|\nabla c|^2}{c} \bigg\} 
	+ \io \frac{|\nabla n|^2}{n} + \io c|D^2 \ln c|^2 \le C \io |u|^4,
	\qquad t>0,
\ee
with some $C>0$ (\cite{win_ctf}).
Appropriate a priori estimates gained from such energy-type inequalities
may allow for the construction of global-in-time solutions, and thereby for going significantly beyond the outcome of the
approach in \cite{lorz}, where local-in-time weak solutions were found for various boundary value problems 
associated with (\ref{0})-(\ref{1}) without making use of such structural properties:
In \cite{DLM}, energy-based arguments were applied to establish global existence of weak solutions, under certain technical
conditions, in the case $\Omega=\R^2$ when $m=1$.
The particular requirement therein that $\|c_0\|_{L^\infty(\R^2)}$ be small was later removed
in \cite{liu_lorz} at the cost of additional structural hypotheses on $\chi$ and $f$.
That also the latter restriction can actually be relaxed was one of the results in \cite{win_ctf},
where a corresponding boundary value problem was considered in bounded convex domains
(cf.~the boundary conditions in (\ref{0b}) below)
under milder assumptions on $\chi$ and $f$, 	
without imposing any smallness condition on the initial data. 
A refined use of energy inequalities was performed there
to derive statements on global existence of weak solutions in the three-dimensional setting and of smooth solutions
in the case $N=2$, even in the situation when the fluid evolution is governed by the full incompressible 
Navier-Stokes equations.
The powerfulness		
of this energy-based approach is further underlined by its ability to moreover yield information
on the large time behavior of solutions if exploited properly: For instance, a generalized version of (\ref{energy}) 
was used in \cite{win_ARMA} as a starting point to show that in the latter two-dimensional chemotaxis-Navier-Stokes 
system all solutions stablilize to the spatially homogeneous equilibrium $(\onz,0,0)$
in the large time limit, where $\onz:=\frac{1}{|\Omega|} \io n_0>0$.\abs
As an apparently inherent drawback, 	
any such type of energy-based reasoning seems to require quite inflexible properties of the parameter functions in (\ref{0}),
thereby possibly excluding even small perturbations of the latter.
For instance, the mentioned results in \cite{liu_lorz} were inter alia built on the strong condition that
$(\frac{f}{\chi})'' < 0$ be valid on $(0,\infty)$.
This could only slightly be relaxed in \cite{win_ctf}, where
merely $(\frac{f}{\chi})'' \le 0$  on $(0,\infty)$ was required; after all, this allowed for the choices
$\chi\equiv const.$ and $f(c)=c$ in (\ref{1}).
Alternative structural conditions can be found in the more recent work \cite{chae_kang_lee}, in which
for the Navier-Stokes variant of (\ref{0}) associated with (\ref{1}) in $\Omega=\R^2$,
a slightly modified version of the Lyapunov functional in (\ref{energy}) is analyzed to
prove global existence of classical solutions 
for nonnegative and noncecreasing $\chi$ and $f$ under the additional condition that
$\|\chi-\mu f\|_{L^\infty((0,\infty))}$ be small for some $\mu\ge 0$.
Even more drastically, for the construction of global weak solutions to the corresponding Cauchy problem in 
$\Omega=\R^3$ it is required there that $\chi$ precisely coincides with a fixed multiple of $f$.\abs
In the case of degenerate cell diffusion of porous medium type, that is, when $m>1$ in (\ref{1}), the above procedure
yields an inequality quite similar to (\ref{energy}), again under essentially the same structural assumptions
on $\chi$ and $f$. Correspondingly obtained a priori estimates can then be exploited
to derive global existence of bounded weak solutions when $N=2$ and $m>1$ is arbitrary (\cite{taowin_DCDS2012}), 
of global weak but possibly unbounded solutions in the case $N=3$ for any $m>1$ (\cite{duan_xiang_IMRN2012}), and of
global weak solutions, locally bounded in $\bar\Omega\times [0,\infty)$, when $N=3$ and $m>\frac{8}{7}$
(\cite{taowin_ANNIHP}), thereby going significantly beyond results achieved without making explicit use thereof
(\cite{DiFLM}, \cite{vorotnikov}).\Abs
{\bf Analysis beyond natural energies.} \quad
More recent experimental findings and corresponding modeling approaches suggest that 
chemotactic migration need not necessarily be directed exclusively toward
increasing signal concentrations, but can rather have rotational components, especially near the physical boundary of 
the domain, and that accordingly the chemotactic sensitivity should actually be considered as a tensor with possibly
nontrivial off-diagonal entries (\cite{xue_othmer}).
In light of the observation that spontaneous emergence of structures indeed seems to occur mainly near droplet 
boundaries (\cite{goldstein2005}),
it thus appears adequate to allow the parameter function $S$ in (\ref{0}) to attain values in $\R^{N\times N}$,
and to thereby depart from the particular structure in (\ref{1}). 
In this general situation, however, it seems no longer possible to derive inequalities of type (\ref{energy}) 
by means of any procedure which in a subtle way cancels contributions stemming from cross-diffusive interaction
as described above.\abs
Accordingly, the goal of the present work will be to develop an alternative a priori estimation method which
is sufficiently robust so as to apply to (\ref{0}) under very mild conditions on all parameter functions appearing therein.
In fact, it turns out that our approach will provide integral estimates which will not only allow for the construction of
global solutions to (\ref{0}) that remain bounded for all times in a suitable sense,
but which beyond this will also serve as a fundament for determining the large time behavior of these solutions.\Abs
{\bf Main results.} \quad
In order to formulate our main results in this direction, let us specify the precise evolution problem addressed in the sequel
by considering (\ref{0}) along with the initial conditions
\be{0i}
    n(x,0)=n_0(x), \quad c(x,0)=c_0(x) \quad \mbox{and} \quad u(x,0)=u_0(x), \qquad x\in\Omega,
\ee
and under the boundary conditions
\be{0b}
    	\Big(D(n)\nabla n - nS(x,n,c)\cdot \nabla c\Big) \cdot \nu=0, \quad \frac{\partial c}{\partial\nu}=0 
	\quad \mbox{and} \quad u=0
	\qquad \mbox{on } \pO,
\ee
in a bounded convex domain $\Omega\subset \R^3$ with smooth boundary,
where throughout this paper we assume for convenience that
\be{init}
    \left\{
    \begin{array}{l}
    n_0 \in C^\kappa(\bar\Omega) \quad \mbox{for some $\kappa>0$ with } n_0\ge 0 \mbox{ in } \Omega, \quad
	\mbox{that} \\
    c_0 \in W^{1,\infty}(\Omega) \quad \mbox{ satisfies $c_0\ge 0$ in $\Omega$, \quad and that}\\
    u_0 \in D(A_r^\beta) \quad \mbox{for some $\beta\in (\frac{1}{2},1)$ and all } r\in (1,\infty).
    \end{array}
    \right.
\ee
with $A_r$ denoting the Stokes operator with domain $D(A_r):=W^{2,r}(\Omega)\cap W_0^{1,r}(\Omega) \cap L^r_\sigma(\Omega)$,
where $L^r_\sigma(\Omega):=\{ \varphi\in L^r(\Omega) \ | \ \nabla \cdot \varphi=0 \}$ for $r\in (1,\infty)$
(cf.~also Section \ref{sect_stokes} below).\\
As for the diffusion coefficient in (\ref{0}), we shall assume that $D$ generalizes the porous-medium-like prototype
$D(n)=mn^{m-1}$ by satisfying 
\be{D1}
    D \in C^{\theta}_{loc}([0,\infty)) \qquad \mbox{for some } \theta>0,
\ee
as well as
\be{D2}
    D(n) \ge \kd n^{m-1} \qquad \mbox{for all } n>0
\ee
with some $m>1$ and $\kd>0$, noting that this includes both degenerate and non-degenerate diffusion at $n=0$.\\
Apart from this, we shall merely suppose that 
\be{S1}
	S\in C^2\Big(\bar\Omega \times [0,\infty)^2; \R^{3\times 3}\Big)
\ee
satisfies
\be{S2}
	|S(x,n,c)| \le S_0(c)
	\quad \mbox{for all $(x,n,c) \in \bar\Omega \times [0,\infty)^2$}
	\qquad \mbox{with some nondecreasing } S_0:[0,\infty)\to \R,
\ee
that
\be{f}
	f\in C^1([0,\infty))
	\quad \mbox{is nonnegative},
\ee
and that
\be{phi}
    \phi \in W^{1,\infty}(\Omega),
\ee
underlining that unlike most previous studies, we do not require any monotonicity property of neither $f$ nor $S$.\abs
In the context of these assumptions, the first of our main results asserts global existence of a bounded solution 
in the following sense.
\begin{theo}\label{theo16}
  Let (\ref{S1}), (\ref{S2}), (\ref{f}) and (\ref{phi}) hold, and suppose that $D$ satisfies
  (\ref{D1}) and (\ref{D2}) with some 
  \be{m}
	m>\frac{7}{6}.
  \ee
  Then for any choice of $n_0, c_0$ and $u_0$ fulfilling (\ref{init}), the problem (\ref{0}), (\ref{0i}), (\ref{0b}) possesses
  at least one global weak solution $(n,c,u,P)$ in the sense of Definition \ref{defi_weak} below.
  This solution is bounded in $\Omega\times (0,\infty)$ in the sense that
  with some $C>0$ we have
  \be{16.1}
    	\|n(\cdot,t)\|_{L^\infty(\Omega)} 
	+ \|c(\cdot,t)\|_{W^{1,\infty}(\Omega)}
	+ \|u(\cdot,t)\|_{W^{1,\infty}(\Omega)}
	\le C \qquad \mbox{for all } t>0.
  \ee
  Furthermore, $c$ and $u$ are continuous in $\bar\Omega\times [0,\infty)$ and
  \be{16.2}
	n \in C^0_{w-\star} ([0,\infty);L^\infty(\Omega));
  \ee
  that is, $n$ is continuous on $[0,\infty)$ as an $L^\infty(\Omega)$-valued function with respect to the 
  weak-$\star$ topology.
\end{theo}
We note here that as compared to the global existence result in \cite{taowin_ANNIHP}, the admissible range for $m$
indicated by (\ref{m}) is slightly smaller. However, besides requiring significantly less conditions on $S$ and $f$,
the statement in Theorem \ref{theo16} goes considerably beyond the outcome in \cite{taowin_ANNIHP} in that, inter alia,
{\em global} boundedness, rather than merely local boundedness, of solutions is obtained here.
A further example for a lack of global boundedness in a system of type (\ref{0})
can be found in \cite{cao_wang}, where global existence of possibly unbounded
classical solutions is proved for the three-dimensional variant  
of (\ref{0}) obtained on specifying $D\equiv 1$ and $f(c)=c$ for $c\ge 0$, and requiring that instead of (\ref{S2}),
$S$ decays for large values of $n$ in the sense that $|S(x,n,c)| \le C(1+n^{-\alpha})$ for all $n\ge 0$ and some $C>0$
and $\alpha>\frac{1}{6}$.\abs
We can moreover show that all the above solutions approach
the unique spatially homogeneous equilibrium corresponding to the bacterial mass $\io n_0$ in the large time limit,
provided that $n_0$ is nontrivial, and that the condition $f>0$ on $(0,\infty)$ is satisfied, 
which is evidently necessary for such a behavior:
\begin{theo}\label{theo166}
  Let (\ref{S1}), (\ref{S2}), (\ref{f}) and (\ref{phi}) hold, suppose that $D$ satisfies
  (\ref{D1}) and (\ref{D2}) with some $m>\frac{7}{6}$, and assume that in addition 
  \be{f_pos}
	f(c)>0
	\qquad \mbox{for all } c>0.
  \ee
  Then whenever $(n_0,c_0,u_0)$ satisfies (\ref{init}) with $n_0\not\equiv 0$,
  the global weak solution constructed in Theorem \ref{theo16} satisfies
  \be{conv}
	n(\cdot,t) \wsto \nzb
	\quad \mbox{ in } L^\infty(\Omega),
	\qquad
	c(\cdot,t) \to 0
	\quad \mbox{ in } L^\infty(\Omega)
	\qquad \mbox{and} \qquad
	u(\cdot,t) \to 0
	\quad \mbox{ in } L^\infty(\Omega)
  \ee
  as $t\to\infty$, where $\nzb:=\frac{1}{|\Omega|} \io n_0$.
\end{theo}
In particular, Theorem \ref{theo16} and Theorem \ref{theo166} provide some progress also in the fluid-free
subcase of (\ref{0}) obtained on letting $\phi\equiv 0$ and $u\equiv 0$. 
Indeed, for the correspondingly gained chemotaxis system with matrix-valued sensitivity the literature so far only
contains very few results: Global classical solutions are known to exist in the case $N=2$ but for small values of
$\|c_0\|_{L^\infty(\Omega)}$ only (\cite{lsxw}); in the same two-dimensional setting, 
global bounded weak solutions for large initial data can be found 
whenever $m>1$ (\cite{cao_ishida_NON}); in the case $m=1$, certain global generalized solutions,
possibly unbounded, have recently been constructed for any $N\ge 1$ and arbitrarily large initial data 
(\cite{win_ct_rot2}). 
If this two-component system is further simplified by moreover choosing the sensitivity to be a constant scalar,
again some energy arguments become available so as to yield much more comprehensive results on global existence of solutions
and even on their large time behavior (\cite{taowin_consumption}).\abs
Our approach underlying the derivation of Theorem \ref{theo16} 	
consists at its core in an analysis of the functional
\be{funct}
	y(t):=\io n^p(\cdot,t) + \io |\nabla c(\cdot,t)|^{2q} + \io |A^\frac{1}{2} u(\cdot,t)|^2,
	\qquad t \ge 0,
\ee
for solutions of certain regularized versions of (\ref{0}) (see Section \ref{sect_approx}), 
where we eventually intend to choose $p>1$ arbitrarily large.
For this purpose, in Section \ref{sect_testing} we shall first follow standard testing procedures to
gain some basic information on the time evolution of each of the summands in (\ref{funct}) separately. 
It will turn out in Section \ref{sect_estimate} and in Section \ref{sect_combine} 
that under suitable conditions on the relationship between the exponents
$p>1$ and $q>1$ herein		
it is possible to estimate the respective ill-signed contributions appropriately, and thereby
establish an ODI for $y$ containing an absorptive linear 
term and thus implying an upper bound for $y$ (see Lemma \ref{lem12} and also Lemma \ref{lem11}).\abs
Since boundedness of $y$ even for large $p$ and $q$ does not directly entail sufficient regularity properties of $u$,
our reasoning will involve a two-step bootstrap argument: In the first step thereof we shall only rely on the 
natural mass conservation property (\ref{mass}) and the smoothing action of the Stokes semigroup
(Section \ref{sect_stokes}) 
to gain a first integral bound for the Jacobian $Du$ (cf.~Lemma \ref{lem13} and also Corollary \ref{cor2}); 
since such bounds allow for estimating certain integrals stemming from the signal-fluid interaction,
provided that $q$ in (\ref{funct}) is not too large (Lemma \ref{lem7} and Lemma \ref{lem8}), this primary information
can be used to derive a bound on $y$ for certain small $p$ and $q$ (Lemma \ref{lem13}). The latter in turn 
implies higher regularity features of $Du$ and thus enables us to treat $y$ for arbitrarily large $p$ and $q$ in a second
step (Lemma \ref{lem14}).\abs
In deriving the convergence properties asserted in Theorem \ref{theo166} in Section \ref{sect_conv}, 
we shall essentially make use of the boundedness statement from Theorem \ref{theo16}: In fact, the latter will
enable us to exploit the finiteness of
$\int_0^\infty \io nf(c)$ and $\int_0^\infty \io |\nabla c|^2$ (Lemma \ref{lem21}) to firstly obtain boundedness also of
$\int_0^\infty \io |\nabla n^\alpha|^2$ for some $\alpha>1$ (Lemma \ref{lem22}), 
and to secondly prove that as a consequence of
these three integral inequalities, all our solutions asymptotically become homogeneous in space and thus satisfy
(\ref{conv}) (Lemma \ref{lem30}, Lemma \ref{lem29} and Lemma \ref{lem31}).
\mysection{Approximation by non-degenerate problems}\label{sect_approx}
Our goal is to construct solutions of (\ref{0}) as limits of solutions to appropriately regularized problems.
To achieve this, we approximate the diffusion coefficient function in (\ref{0}) by introducing a family
$(D_\eps)_{\eps\in (0,1)}$ of functions
\bas
	& & D_\eps \in C^2([0,\infty))
	\quad \mbox{such that} \quad
	D_\eps(n) \ge \eps \quad \mbox{for all $n\ge 0$ and } \eps\in (0,1)
	\quad \mbox{and} \quad \\
	& & D(n) \le D_\eps(n) \le D(n)+2\eps
	\quad \mbox{for all $n\ge 0$ and } \eps\in (0,1).
\eas
Next, it will be convenient to deal with homogeneous Neumann boundary conditions for both $n$ and $c$
rather than with the nonlinear no-flux relation in (\ref{0b}). 
In order to achieve this at least during our approximation procedure, following \cite{lsxw} we moreover
fix families $(\rho_\eps)_{\eps\in (0,1)}$ and $(\chi_\eps)_{\eps\in (0,1)}$ of functions 
\bas
	\rho_\eps \in C_0^\infty(\Omega)
	\quad \mbox{with} \quad
	0 \le \rho_\eps \le 1 \mbox{ in $\Omega$ \quad and \quad}
	\rho_\eps \nearrow 1 \mbox{ in $\Omega$ as $\eps\searrow 0$}
\eas
and
\bas
	\chi_\eps \in C_0^\infty([0,\infty))
	\quad \mbox{satisfying} \quad
	0 \le \chi_\eps \le 1 \mbox{ in $[0,\infty)$ \quad and \quad}
	\chi_\eps \nearrow 1 \mbox{ in $[0,\infty)$ as $\eps\searrow 0$,}
\eas
and define smooth approximations $S_\eps \in C^2(\bar\Omega\times [0,\infty)^2;\R^{3\times 3})$ of $S$ by letting
\be{Seps}
	S_\eps(x,n,c):=\rho_\eps(x) \cdot \chi_\eps(n) \cdot S(x,n,c),
	\qquad x\in\bar\Omega, \ n\ge 0, \ c \ge 0,
\ee
for $\eps\in (0,1)$.
Then for any such $\eps$, the regularized problems
\be{0eps}
	\left\{ \begin{array}{ll}
	n_{\eps t} + \ueps \cdot \nabla \neps = \nabla \cdot \Big(D_\eps(\neps)\nabla \neps \Big) 
	- \nabla \cdot \Big(\neps S_{\eps}(x,\neps,\ceps) \cdot \nabla \ceps \Big),
	\qquad & x\in \Omega, \ t>0, \\[1mm]
	c_{\eps t}+ \ueps \cdot \nabla \ceps = \Delta \ceps - \neps f(\ceps),
	\qquad & x\in \Omega, \ t>0, \\[1mm]
	u_{\eps t} + \nabla \Peps = \Delta \ueps + \neps \nabla \phi,
	\qquad & x\in \Omega, \ t>0, \\[1mm]
	\nabla \cdot \ueps=0,
	\qquad & x\in \Omega, \ t>0, \\[1mm]
	\frac{\partial \neps}{\partial\nu} = \frac{\partial \ceps}{\partial\nu}=0, \quad \ueps=0,
	\qquad & x\in \partial\Omega, \ t>0, \\[1mm]
	\neps(x,0)=n_0(x), \quad \ceps(x,0)=c_0(x), \quad \ueps(x,0)=u_0(x),
	\qquad & x\in\Omega,
 	\end{array} \right.
\ee
are globally solvable in the classical sense:
\begin{lem}\label{lem0}
  Let $\eps\in (0,1)$. Then there exist functions
  \bas
	\left\{ \begin{array}{l}
	\neps\in C^0(\bar\Omega \times [0,\infty)) \cap C^{2,1}(\bar\Omega\times (0,\infty)), \\
	\ceps\in C^0(\bar\Omega \times [0,\infty)) \cap C^{2,1}(\bar\Omega\times (0,\infty)), \\
	\ueps\in C^0(\bar\Omega \times [0,\infty)) \cap C^{2,1}(\bar\Omega\times (0,\infty)), \\
	\Peps\in C^{1,0}(\bar\Omega\times (0,\infty)),
	\end{array} \right.
  \eas
  such that $(\neps,\ceps,\ueps,\Peps)$ solves (\ref{0eps}) classically in $\Omega \times (0,\infty)$, and such
  that $\neps$ and $\ceps$ are nonnegative in $\Omega\times (0,\infty)$.
\end{lem}
\proof
  By an adaptation of well-established fixed point arguments, one can readily verify the existence of a local-in-time
  smooth solution, nonnegative in its first two components by the maximum principle, 
  and extensible up to a maximal time $\tme \in (0,\infty]$ which in the case $\tme<\infty$ 
  has the property that
  \be{ext_eps}
	\limsup_{t\nearrow \tme} \Big( \|\neps(\cdot,t)\|_{C^2(\bar\Omega)}
	+ \|\ceps(\cdot,t)\|_{C^2(\bar\Omega)}
	+ \|\ueps(\cdot,t)\|_{C^2(\bar\Omega)}
	\Big) = \infty
  \ee
  (cf.~\cite[Lemma 2.1]{win_ctf} and \cite[Lemma 2.1]{taowin1}, for instance). 
  Since (\ref{Seps}) ensures that for fixed $\eps\in (0,1)$ the function $S_\eps(x,n,c)$ vanishes for all 
  sufficiently large $n$, one may apply standard a priori estimation techniques to infer that
  for any such $\eps$ and each $T>1$ there exists $C(\eps,T)>0$ such that
  \bas
	\|\neps(\cdot,t)\|_{C^2(\bar\Omega)}
	+ \|\ceps(\cdot,t)\|_{C^2(\bar\Omega)}
	+ \|\ueps(\cdot,t)\|_{C^2(\bar\Omega)}
	\le C(\eps,T)
	\qquad \mbox{for all } t\in (\tau,\tilde \tme),
  \eas
  where $\tau:=\min\{1,\frac{1}{2} \tme\}$ and $\tilde \tme:=\min\{T,\tme\}$
  (cf.~e.g.~\cite[Sect. 5]{win_ctf}, \cite{horstmann_win} and \cite{taowin_JDE}).
  As a consequence of this and (\ref{ext_eps}), we actually must have $\tme=\infty$, as desired.
\qed
The following basic properties of solutions to (\ref{0eps}) are immediate.
\begin{lem}\label{lem_basic}
  The solution of (\ref{0eps}) satisfies
  \be{mass}
	\|\neps(\cdot,t)\|_{L^1(\Omega)}=\io n_0 
	\qquad \mbox{for all } t>0
  \ee
  as well as
  \be{cinfty}
	\|\ceps(\cdot,t)\|_{L^\infty(\Omega)} \le \|c_0\|_{L^\infty(\Omega)}
	\qquad \mbox{for all } t>0.
  \ee
\end{lem}
\proof
  The mass conservation property (\ref{mass}) directly follows by integrating the first equation in (\ref{0eps})
  over $\Omega$. 
  Moreover, using that both $\neps$ and $f$ are nonnegative we can readily derive 
  the inequality (\ref{cinfty}) by applying a parabolic comparison argument to the second equation in (\ref{0eps}).
\qed
\mysection{A priori estimates}
We proceed to derive $\eps$-independent estimates for the approximate solutions constructed above. Throughout this
section, for $\eps\in (0,1)$ we let $(\neps,\ceps,\ueps,P_\eps)$ denote the global solution of (\ref{0eps}).
\subsection{$W^{1,r}$ regularity of $u$ implied by $L^p$ regularity of $n$}\label{sect_stokes}
Our first goal is to draw a consequence of a supposedly known bound for $\neps$ in $L^\infty((0,\infty);L^p(\Omega))$
on the regularity of the spatial derivative $D\ueps$ of $\ueps$. 
Since in view of (\ref{mass}) we intend to apply this inter alia to the case $p=1$
in a first step (cf.~Lemma \ref{lem13}), applying smoothing estimates for the Stokes semigroup seems not fully 
straightforward in our situation. \\
In order to prepare our results in this direction, let us recall that for each $r\in (1,\infty)$, the Helmholtz projection
acts as a bounded linear operator $\proj_r$ from $L^r(\Omega)$ onto its subspace
$L^r_\sigma(\Omega)=\{\varphi\in L^r(\Omega) \ | \ \nabla \cdot \varphi=0\}$ of all solenoidal vector fields.
Moreover, the 
realization $A_r$ of the Stokes operator $A$ in $L^r_\sigma(\Omega)$
with domain $D(A_r)=W^{2,r}(\Omega)\cap W_0^{1,r}(\Omega) \cap L^r_\sigma(\Omega)$ is sectorial in $L^r_\sigma(\Omega)$
and hence possesses closed fractional powers $A_r^\beta$ with dense domains $D(A_r^\beta)$ for any $\beta\in \R$ 
(\cite{giga1981_the_other}),
and $A_r$ generates an analytic semigroup $(e^{-t A_r})_{t\ge 0}$ in $L^r_\sigma(\Omega)$.
In the sequel, 
since $\proj_r \psi$ and $A_r^\beta \varphi$ as well as $e^{-t A_r} \varphi$ are actually independent of $r\in (1,\infty)$ 
for each $\psi\in C_0^\infty(\Omega)$ and 
$\varphi\in C_0^\infty(\Omega) \cap L^r_\sigma(\Omega)$, $\beta\in\R$ and $t\ge 0$, 
we may suppress the subscript $r$ in $\proj_r, A_r^\beta$ and $e^{-t A_r}$ whenever there is no danger of confusion.\\
Then among well-known embedding and regularity estimates we will especially need the following in the sequel.
\begin{lem}\label{lem1}
  i) \ Let $r>1$. Then for all $\beta>\frac{1}{2}$ one can find $C=C(r,\beta)>0$ such that
  \be{1.1}
	\|D\varphi\|_{L^r(\Omega)} \le C \|A^\beta \varphi\|_{L^r(\Omega)}
	\qquad \mbox{for all } \varphi\in D(A_r^\beta).
  \ee
  ii) \ There exists $\mu>0$ with the following property:
  For all $r\in (1,\infty)$ and $p\in (1,r]$ and each $\beta\ge 0$ there exists $C=C(p,r,\beta)$ such that
  whenever $\varphi \in L^p_\sigma(\Omega)$, we have
  \be{1.2}
	\|A^\beta e^{-tA} \varphi\||_{L^r(\Omega)}
	\le C t^{-\beta - \frac{3}{2}(\frac{1}{p}-\frac{1}{r})} e^{-\mu t} \|\varphi\|_{L^p(\Omega)}
	\qquad \mbox{for all } t>0.
  \ee
\end{lem}
\proof
  i) \ The inequality (\ref{1.1}) can be derived from \cite[Lemma II.17.1]{friedman} by adapting the argument in
  \cite[Lemma II.17.2]{friedman} in a straightforward manner.\\
  ii) \ For (\ref{1.2}), we may refer the reader to \cite[p.201]{giga1986}.
\qed
Evidently, the regularity properties for $\ueps$ are linked to those of the forcing term	
$g_\eps:=\proj[\neps\nabla\phi]$
appearing in the version $u_{\eps t}+A \ueps=\proj[\neps\nabla\phi]$ of the Stokes subsystem of 
(\ref{0eps}) when projected to the respective spaces of divergence-free functions.
Since in a first step, our basic information (\ref{mass}) only asserts a bound for $g_\eps$ with respect
to the norm in $L^1(\Omega)$, and since standard results apparently do not apply directly to this 
non-reflexive situation, we briefly include the following two lemmata 
to prepare an adequate estimation of $g_\eps$.\\
Firstly, a further embedding property of the domains of fractional powers of $A$ allows us to control, 
upon a slight lifting through a negative fractional power of $A$, the norm of functions
in the space $L^{p_0}(\Omega)$ by the norm of the unlifted function in $L^p(\Omega)$ with some $p$ smaller than 
$p_0$. 
In Lemma \ref{lem01}, 
the case $p_0=\infty$ will be of particular importance for treating the $L^1$ situation mentioned above.
\begin{lem}\label{lem02}
  Suppose that $1<p< p_0 \le \infty$, and that $\delta\in (0,1)$ is such that $2\delta-\frac{3}{p}>-\frac{3}{p_0}$.
  Then there exists $C>0$ such that
  \be{02.1}
	\|A^{-\delta} \psi\|_{L^{p_0}(\Omega)} \le C \|\psi\|_{L^{p}(\Omega)}
	\qquad \mbox{for all } \psi\in L^{p}(\Omega).
  \ee
\end{lem}
\proof
  According to \cite[Theorem 3]{giga1981_the_other} and
  \cite[Theorem 1.6.1]{henry}, our assumption on $\delta$ ensures that $D(A_p^\delta) \hra L^{p_0}(\Omega)$,
  which means that there exists $C_1>0$ such that
  \bas
	\|\varphi\|_{L^{p_0}(\Omega)} \le C_1 \|A^\delta \varphi\|_{L^p(\Omega)}
	\qquad \mbox{for all } \varphi \in D(A_p^\delta).
  \eas
  Thus, if we fix $\psi\in C_0^\infty(\Omega)$ and apply this to $\varphi:=A^{-\delta}\psi$, we see that the inequality
  in (\ref{02.1}) holds with $C:=C_1$. For arbitrary $\psi\in L^p(\Omega)$, (\ref{02.1}) easily follows from this
  by completion.
\qed
By means of a straightforward duality argument, we can thereby indeed use a knowledge on the size of a function 
in $L^1(\Omega)$ to control a slightly lifted variant of its solenoidal part in a reflexive $L^p$ space.
More generally, we have the following.
\begin{lem}\label{lem01}
  Assume that $1\le p < p_0<\infty$, and that $\delta\in (0,1)$ is such that $2\delta-\frac{3}{p}>-\frac{3}{p_0}$.
  Then there exists $C>0$ such that
  \be{01.1}
	\|A^{-\delta} \proj \psi\|_{L^{p_0}(\Omega)} \le C \|\psi\|_{L^p(\Omega)}
	\qquad \mbox{for all } \psi\in C_0^\infty(\Omega).
  \ee
  Consequently, the operator $A^{-\delta} \proj$ possesses a unique extension to all of $L^{p_0}(\Omega)$
  with norm controlled according to (\ref{01.1}).
\end{lem}
\proof
  Let $\varphi\in C_0^\infty(\Omega)$. Then since both $A^{-\delta} \proj \psi$ and $A^{-\delta} \proj \varphi$ are  
  divergence free and $A^{-\delta}$ is symmetric, we have
  \bas
	\io A^{-\delta}\proj \psi \cdot \varphi
	= \io A^{-\delta} \proj \psi \cdot \proj \varphi 
	= \io \proj \psi \cdot A^{-\delta} \proj \varphi 
	= \io \psi \cdot A^{-\delta} \proj \varphi.
  \eas
  Now with $p':=\frac{p}{p-1} \in (1,\infty]$ and $p_0':=\frac{p_0}{p_0-1}$ we have $1<p_0'<p' \le \infty$,
  and the assumption $2\delta-\frac{3}{p}>-\frac{3}{p_0}$
  ensures that $2\delta-\frac{3}{p_0'}>-\frac{3}{p'}$. Therefore we may invoke Lemma \ref{lem02}
  and use the boundedness of the projection $\proj$ in $L^{p_0'}(\Omega)$ to find $C_1>0$ and $C_2>0$ such that
  \bas
	\bigg| \io A^{-\delta} \proj \psi \cdot \varphi \bigg|
	&\le& \|\psi\|_{L^p(\Omega)} \cdot \|A^{-\delta} \proj \varphi\|_{L^{p'}(\Omega)} \\
	&\le& C_1 \|\psi\|_{L^p(\Omega)} \cdot \|\proj \varphi\|_{L^{p_0'}(\Omega)} \\
	&\le& C_2 \|\psi\|_{L^p(\Omega)} \cdot \|\varphi\|_{L^{p_0'}(\Omega)}
	\qquad \mbox{for all } \varphi\in C_0^\infty(\Omega).
  \eas
  By a standard duality argument, this implies (\ref{01.1}).
\qed
An application of this to the Stokes equations in (\ref{0eps}) now yields the following implication
of some presupposed boundedness property of $\neps$ to the regularity features of $\ueps$.
\begin{cor}\label{cor2}
  Let $p \in [1,\infty)$ and $r\in [1,\infty]$ be such that
  \be{2.1}
	\left\{ \begin{array}{ll}
	r < \frac{3p}{3-p} \qquad & \mbox{if } p\le 3, \\[1mm]
	r \le \infty \qquad & \mbox{if } p>3.
 	\end{array} \right.
  \ee
  Then for all $K>0$ there exists $C=C(p,r,K)$ such that if for some $\eps\in (0,1)$ and $T>0$ we have
  \be{2.2}
	\|\neps(\cdot,t)\|_{L^p(\Omega)} \le K
	\qquad \mbox{for all } t\in (0,T),
  \ee
  then
  \be{2.3}
	\|D\ueps(\cdot,t)\|_{L^r(\Omega)} \le C
	\qquad \mbox{for all } t\in (0,T).
  \ee
\end{cor}
\proof
  In view of (\ref{2.1}), it is evidently sufficient to consider the case $r>p$ only, in which we can fix $r_0\in (p,r)$
  such that with $\beta$ as in (\ref{init}) we have
  \be{2.33}
	\frac{1}{2} + \frac{3}{2}\Big( \frac{1}{r_0}-\frac{1}{r}\Big) < \beta.
  \ee
  Since (\ref{2.1}) moreover ensures that
  \bas
	\bigg\{ \frac{1}{2} + \frac{3}{2}\Big(\frac{1}{r_0}-\frac{1}{r}\Big) \bigg\}
	- \bigg\{ 1 - \frac{3}{2}\Big(\frac{1}{p}-\frac{1}{r_0}\Big)\bigg\}
	= -\frac{1}{2}+\frac{3}{2}\Big(\frac{1}{p}-\frac{1}{r}\Big)<0,
  \eas
  we can thus choose $\beta_0\in (\frac{1}{2},\beta)$ fulfilling
  \be{2.4}
	\frac{1}{2} + \frac{3}{2}\Big( \frac{1}{r_0}-\frac{1}{r}\Big) < \beta_0
	< 1 - \frac{3}{2}\Big(\frac{1}{p}-\frac{1}{r_0}\Big),
  \ee
  then pick $\delta \in (0,1)$ small enough such that still
  \be{2.44}
	\beta_0+\delta	
	< 1 - \frac{3}{2}\Big(\frac{1}{p}-\frac{1}{r_0}\Big),
  \ee
  and finally fix some $p_0>p$ sufficiently close to $p$ such that 
  \be{2.444}
	2\delta-\frac{3}{p}>-\frac{3}{p_0}.
  \ee
  Then in the variation-of-constants representation
  \be{2.45}
	\ueps(\cdot,t) = e^{-tA} u_0
	+ \int_0^t e^{-(t-s)A} \proj [\neps(\cdot,s)\nabla \phi] ds,
	\qquad t\in (0,T),
  \ee
  with $A=A_{r_0}$ and $\proj=\proj_{r_0}$, we apply $A^{\beta_0}$ on both sides and use Lemma \ref{lem1} i) along with the
  left inequality in (\ref{2.4}) to find $C_1>0$ such that
  \bea{2.5}
	\hspace*{-10mm}
	\|D\ueps(\cdot,t)\|_{L^r(\Omega)}
	&\le& C_1 \|A^{\beta_0} \ueps(\cdot,t)\|_{L^{r_0}(\Omega)} \nn\\
	&\le& C_1 \|A^{\beta_0} e^{-tA} u_0\|_{L^{r_0}(\Omega)}
	+ C_1 \int_0^t \Big\| A^{\beta_0+\delta} e^{-(t-s)A} 
	A^{-\delta} \proj [\neps(\cdot,s) \nabla \phi] \Big\|_{L^{r_0}(\Omega)} ds
  \eea
  for all $t\in (0,T)$. Here since $u_0\in D(A^{\beta_0})$ by (\ref{init}) and the fact that $\beta_0<\beta$, we have
  \bas
	\|A^{\beta_0} e^{-tA} u_0\|_{L^{r_0}(\Omega)}
	= \|e^{-tA} A^{\beta_0} u_0\|_{L^{r_0}(\Omega)}
	\le C_2
	\qquad \mbox{for all } t\in (0,T)
  \eas
  with some $C_2>0$.
  Furthermore, since $p_0>p\ge 1$, Lemma \ref{lem1} ii) applies to show that there exist $C_3>0$ and $\mu>0$ fulfilling

  \bas
	\Big\| A^{\beta_0+\delta} e^{-(t-s)A} 
	A^{-\delta} \proj [\neps(\cdot,s) \nabla \phi] \Big\|_{L^{r_0}(\Omega)}
	\le C_3(t-s)^{-\beta_0-\delta-\frac{3}{2}(\frac{1}{p_0}-\frac{1}{r_0})} e^{-\mu(t-s)}
	\Big\|A^{-\delta} \proj[\neps(\cdot,s)\nabla \phi]\Big\|_{L^{p_0}(\Omega)}
  \eas
  for all $t\in (0,T)$ and $s\in (0,t)$,
  where thanks to Lemma \ref{lem01} and the boundedness of $\nabla \phi$, we can use (\ref{2.2}) to see that
  \bas
	\Big\|A^{-\delta} \proj[\neps(\cdot,s)\nabla \phi]\Big\|_{L^{p_0}(\Omega)}
	\le C_4 \|\neps(\cdot,s)\nabla\phi\|_{L^p(\Omega)}
	\le C_5 \|\neps(\cdot,s)\|_{L^p(\Omega)}
	\le C_5 K
	\qquad \mbox{for all } s\in (0,T)
  \eas
  with some $C_4>0$ and $C_5>0$.
  Therefore, from (\ref{2.5}) we altogether obtain that
  \bas
	\|D\ueps(\cdot,t)\|_{L^r(\Omega)}
	&\le& C_1 C_2 + C_1 C_3 C_5 K \int_0^t (t-s)^{-\beta_0-\delta-\frac{3}{2}(\frac{1}{p_0}-\frac{1}{r_0})}
	e^{-\mu(t-s)} ds \\
	&\le& C_1 C_2 + C_1 C_3 C_5 C_6 K
	\qquad \mbox{for all } t\in (0,T),
  \eas
  where
  \bas
	C_6 := \int_0^\infty \sigma^{-\beta_0-\delta-\frac{3}{2}(\frac{1}{p_0}-\frac{1}{r_0})} e^{-\mu\sigma} d\sigma
  \eas
  is finite according to (\ref{2.44}). This proves (\ref{2.3}).
\qed
\subsection{Standard testing procedures}\label{sect_testing}
We now turn to the analysis of the coupled functional in (\ref{funct}).
Here we first apply standard testing procedures to gain the inequalities in the following three lemmata.
Further estimating the respective right-hand sides therein will then be done separately in the sequel.\\
Let us begin by testing the first equation in (\ref{0eps}) against powers of $\neps$.
\begin{lem}\label{lem4}
  Let $p>1$. Then for all $\eps\in (0,1)$,
  \be{4.1}
	\frac{1}{p} \frac{d}{dt} \io \neps^p 
	+ \frac{2(p-1)\kd}{(m+p-1)^2} \io |\nabla \neps^\frac{p+m-1}{2}|^2
	\le \frac{(p-1)S_1^2}{2\kd} \io \neps^{p+1-m} |\nabla \ceps|^2
	\qquad \mbox{for all } t>0,
  \ee
  where $\kd$ is as in (\ref{D2}) and
  \bas
	S_1:=S_0(\|c_0\|_{L^\infty(\Omega)})
  \eas
  with $S_0$ taken from (\ref{S2}).
\end{lem}
\proof
  We multiply the first equation in (\ref{0eps}) by $\neps^{p-1}$ and integrate by parts over $\Omega$. 
  Since $S_\eps(x,\neps,\ceps)$
  vanishes whenever $x\in\pO$ according to (\ref{Seps}), this yields
  \be{4.2}
	\frac{1}{p} \frac{d}{dt} \io \neps^p + (p-1) \io \neps^{p-2} D_\eps(\neps) |\nabla \neps|^2
	= (p-1) \io \neps^{p-1} \nabla \neps \cdot \Big(S_\eps(x,\neps,\ceps) \cdot \nabla \ceps\Big)
	\qquad \mbox{for all } t>0.
  \ee
  Here we use the definition of $D_\eps$ and (\ref{D2}) to see that
  \bas
	(p-1)\io \neps^{p-2} D_\eps(\neps) |\nabla \neps|^2
	\ge (p-1) \kd \io \neps^{p+m-3} |\nabla \neps|^2
	\qquad \mbox{for all } t>0,
  \eas  
  and next combine (\ref{Seps}) with (\ref{S2}) and (\ref{cinfty}) to obtain
  \bas
	|S_\eps(x,\neps,\ceps)| \le S_1
	\qquad \mbox{in } \Omega\times (0,\infty),
  \eas
  so that using (\ref{D2}) and Young's inequality we can estimate
  \bas
	(p-1) \io \neps^{p-1} \nabla \neps \cdot \Big(S_\eps(x,\neps,\ceps)\cdot \nabla \ceps \Big)
	&\le& (p-1)S_1 \io \neps^{p-1} |\nabla \neps| \cdot |\nabla \ceps| \\
	&\le& \frac{(p-1)\kd}{2} \io \neps^{p+m-3} |\nabla \neps|^2
	+ \frac{(p-1) S_1^2}{2\kd} \io \neps^{p+1-m} |\nabla \ceps|^2
  \eas
  for all $t>0$, whence (\ref{4.1}) readily follows from (\ref{4.2}).
\qed
In order to obtain a first information on the time evolution also of 
$\io |\nabla \ceps|^{2q}$, we make use of the convexity of $\Omega$ in the following lemma.
\begin{lem}\label{lem5}
  Let $q>1$ and $\eps\in (0,1)$. Then
  \bea{5.01}
	& & \hspace*{-20mm}
	\frac{1}{2q} \frac{d}{dt} \io |\nabla \ceps|^{2q}
	+ \frac{2(q-1)}{q^2} \io \Big| \nabla |\nabla \ceps|^q \Big|^2 
	+ \frac{1}{2} \io |\nabla \ceps|^{2q-2} |D^2 \ceps|^2 \nn\\
	&\le& \frac{(2q-2+\sqrt{3})^2 f_1}{2} \io \neps^2 |\nabla \ceps|^{2q-2}
	+ \io |\nabla \ceps|^{2q} \cdot |D\ueps|
	\qquad \mbox{for all } t>0,
  \eea
  where 
  \be{f1}
	f_1:=\|f\|_{L^\infty((0,\|c_0\|_{L^\infty(\Omega)}))}.
  \ee
\end{lem}
\proof
  Differentiating the second equation in (\ref{0eps}) and using that 
  $\Delta |\nabla \ceps|^2=2\nabla \ceps\cdot\nabla \Delta \ceps + 2|D^2 \ceps|^2$, we obtain the pointwise identity
  \bea{5.1}
	\hspace*{-10mm}
	\frac{1}{2} \Big(|\nabla \ceps|^2\Big)_t
	&=& \nabla \ceps \cdot \nabla \Big\{ \Delta \ceps - \neps f(\ceps)- \ueps\cdot \nabla \ceps\Big\} \nn\\
	&=& \frac{1}{2} \Delta |\nabla \ceps|^2 - |D^2 \ceps|^2
	- \nabla \ceps \cdot \nabla (\neps f(\ceps))
	- \nabla \ceps \cdot \nabla (\ueps\cdot\nabla \ceps)
	\quad \mbox{in } \Omega\times (0,\infty).
  \eea
  We multiply this by $(|\nabla \ceps|^2)^{q-1}$ and integrate by parts over $\Omega$.
  Since $\frac{\partial \ceps}{\partial \nu}=0$ on $\pO$ along with the convexity of $\Omega$ ensures that
  $\frac{\partial |\nabla \ceps|^2}{\partial\nu} \le 0$ on $\pO$ (\cite[Lemma I.1]{lions_ARMA}), 
  this results in the inequality
  \bea{5.2}
	& & \hspace*{-20mm}
	\frac{1}{2q} \frac{d}{dt} \io |\nabla \ceps|^{2q}
	+ \frac{q-1}{2} \io |\nabla \ceps|^{2q-4} \Big|\nabla |\nabla \ceps|^2 \Big|^2
	+ \io |\nabla \ceps|^{2q-2} |D^2 \ceps|^2 \nn\\
	&\le& - \io |\nabla \ceps|^{2q-2} \nabla \ceps \cdot \nabla (\neps f(\ceps))
	- \io |\nabla \ceps|^{2q-2} \nabla \ceps \cdot \nabla (\ueps\cdot \nabla \ceps)
	\qquad \mbox{for all } t>0.
  \eea
  Here in the first integral on the right we again integrate by parts to estimate, using that 
  \bas
	|f(\ceps)| \le f_1
	\qquad \mbox{in } \Omega\times (0,\infty)
  \eas 
  by (\ref{cinfty}) and (\ref{f1}),
  \bas
	- \io |\nabla \ceps|^{2q-2} \nabla \ceps \cdot \nabla (\neps f(\ceps))
	&=& \io \neps f(\ceps) |\nabla \ceps|^{2q-2} \Delta \ceps
	+ \io \neps f(\ceps) \nabla \ceps \cdot \nabla |\nabla \ceps|^{2q-2} \\
	&\le& f_1 \io \neps |\nabla \ceps|^{2q-2} |\Delta \ceps|
	+ f_1 \io \neps |\nabla \ceps| \cdot \Big| \nabla |\nabla \ceps|^{2q-2} \Big|
  \eas
  for all $t>0$. Since $|\Delta \ceps| \le \sqrt{3} |D^2 \ceps|$ by the Cauchy-Schwarz inequality, and since
  \bas
	\nabla |\nabla \ceps|^{2q-2}
	= 2(q-1) |\nabla \ceps|^{2q-4} D^2 \ceps \cdot \nabla \ceps
	\qquad \mbox{in } \Omega\times (0,\infty),
  \eas
  in view of Young's inequality this implies that
  \bea{5.3}
	- \io |\nabla \ceps|^{2q-2} \nabla \ceps \cdot \nabla (\neps f(\ceps))
	&\le& \sqrt{3} f_1 \io \neps |\nabla \ceps|^{2q-2} |D^2 \ceps|
	+ 2(q-1) f_1 \io \neps |\nabla \ceps|^{2q-3} |D^2 \ceps \cdot \nabla \ceps| \nn\\
	&\le& (2q-2+\sqrt{3}) f_1 \io \neps |\nabla \ceps|^{2q-2} |D^2 \ceps| \nn\\
	&\le& \frac{1}{2} \io |\nabla \ceps|^{2q-2} |D^2 \ceps|^2
	+ \frac{(2q-2+\sqrt{3})^2 f_1^2}{2} \io \neps^2 |\nabla \ceps|^{2q-2}
  \eea
  for all $t>0$.
  As for the rightmost integral in (\ref{5.2}), we first differentiate $\ueps\cdot \nabla \ceps$ 
  to gain the decomposition
  \bea{5.4}
	\hspace*{-5mm}
	- \io |\nabla \ceps|^{2q-2} \nabla \ceps \cdot \nabla (\ueps\cdot \nabla \ceps)
	=- \io |\nabla \ceps|^{2q-2} \nabla \ceps \cdot (D\ueps \cdot \nabla \ceps)
	- \io |\nabla \ceps|^{2q-2} \nabla \ceps \cdot (D^2 \ceps \cdot \ueps)
  \eea
  for all $t>0$, and then use the pointwise identity
  \bas
	|\nabla \ceps|^{2q-2} \nabla \ceps \cdot (D^2 \ceps \cdot \ueps)
	&=& \ueps \cdot \Big( |\nabla \ceps|^{2q-2} D^2 \ceps \cdot \nabla \ceps\Big) \\
	&=& \frac{1}{2q} \, \ueps\cdot \nabla |\nabla \ceps|^{2q}
	\qquad \mbox{in } \Omega\times (0,\infty)
  \eas
  to infer on integrating by parts that
  \bas
	-\io |\nabla \ceps|^{2q-2} \nabla \ceps \cdot (D^2 \ceps \cdot \ueps)
	&=& - \frac{1}{2q} \io \ueps\cdot \nabla |\nabla \ceps|^{2q} \\
	&=& \frac{1}{2q} \io (\nabla \cdot \ueps) |\nabla \ceps|^{2q} \\[2mm]
	&=& 0
	\qquad \mbox{for all } t>0,
  \eas
  because $\ueps$ is solenoidal. Therefore, (\ref{5.4}) entails that
  \bas
	- \io |\nabla \ceps|^{2q-2} \nabla \ceps \cdot \nabla (\ueps\cdot \nabla \ceps)
	&=& - \io |\nabla \ceps|^{2q-2} \nabla \ceps \cdot (D\ueps \cdot \nabla \ceps) \\
	&\le& \io |\nabla \ceps|^{2q} \cdot |D\ueps|
	\qquad \mbox{for all } t>0,
  \eas
  which combined with (\ref{5.3}) and (\ref{5.2}) yields (\ref{5.01}).
\qed
Finally, the following basic inequality for $\io |A^\frac{1}{2}\ueps|^2$ is standard.
\begin{lem}\label{lem55}
  For any $\eps\in (0,1)$, we have
  \be{55.1}
	\frac{d}{dt} \io |A^\frac{1}{2} \ueps|^2
	+ \io |A\ueps|^2
	\le \|\nabla \phi\|_{L^\infty(\Omega)}^2 \cdot \io \neps^2
	\qquad \mbox{for all } t>0.
  \ee
\end{lem}
\proof

  We apply the Helmholtz projection $\proj$ to the third equation in (\ref{0eps}) and test the resulting identity
  $u_{\eps t}+A\ueps = \proj(\neps\nabla\phi)$ by $A\ueps$.
  Using that $A^\frac{1}{2}$ is self-adjoint in $L^2_\sigma(\Omega)$ and that $\proj$ acts as  
  an orthogonal projection in this Hilbert space,
  by Young's inequality we obtain that
  \bas
	\frac{1}{2} \frac{d}{dt} \io |A^\frac{1}{2} \ueps|^2 + \io |A\ueps|^2
	&=& \io A\ueps \cdot \proj (\neps\nabla\phi) \\
	&\le& \frac{1}{2} \io |A\ueps|^2 + \frac{1}{2} \io |\proj (\neps\nabla\phi)|^2 \\
	&\le& \frac{1}{2} \io |A\ueps|^2 + \frac{1}{2} \io |\neps \nabla \phi|^2 \\
	&\le& \frac{1}{2} \io |A\ueps|^2 + \frac{1}{2} \|\nabla\phi\|_{L^\infty(\Omega)}^2 \io \neps^2
	\qquad \mbox{for all } t>0,
  \eas
  and thereby precisely arrive at (\ref{55.1}).
\qed
\subsection{Estimating the right-hand sides in (\ref{4.1}), (\ref{5.01}) and (\ref{55.1})}\label{sect_estimate}
We next plan to estimate the right-hand sides in the above inequalities appropriately
by using suitable interpolation arguments along with the basic a priori information provided by Lemma \ref{lem_basic}.
Here the following auxiliary interpolation lemma, extending a similar statement known in the 
two-dimensional case (\cite{lsxw}) to the present framework, 
will play an important role in making efficient use of the known $L^\infty$ bound (\ref{cinfty}) for $\ceps$.
\begin{lem}\label{lem6}
  Suppose that $\Omega\subset \R^3$ is a bounded domain with smooth boundary, that $q\ge 1$ and that
  \be{6.01}
	\lambda \in [2q+2,4q+1].
  \ee
  Then there exists $C>0$ such that for all $\varphi\in C^2(\bar\Omega)$ fulfilling 
  $\varphi\cdot \frac{\partial\varphi}{\partial\nu}=0$ on $\pO$ we have
  \be{6.1}
	\|\nabla \varphi\|_{L^\lambda(\Omega)}
	\le C \Big\| |\nabla \varphi|^{q-1} D^2\varphi \Big\|_{L^2(\Omega)}^\frac{2\lambda-6}{(2q-1)\lambda}
	\|\varphi\|_{L^\infty(\Omega)}^\frac{6q-\lambda}{(2q-1)\lambda} 
	+ C \|\varphi\|_{L^\infty(\Omega)}.
  \ee
\end{lem}
\proof
  We integrate by parts to rewrite
  \bea{6.2}
	\|\nabla\varphi\|_{L^\lambda(\Omega)}^\lambda
	= - \io \varphi |\nabla \varphi|^{\lambda-2} \Delta\varphi
	- (\lambda-2) \io \varphi |\nabla\varphi|^{\lambda-4} \nabla\varphi \cdot (D^2 \varphi \cdot\nabla\varphi),
  \eea
  where by the Cauchy-Schwarz inequality we see that
  \be{6.3}
	\bigg|	- (\lambda-2) \io \varphi |\nabla\varphi|^{\lambda-4} \nabla\varphi \cdot (D^2 \varphi \cdot\nabla\varphi)
	\bigg|
	\le (\lambda-2) \|\varphi\|_{L^\infty(\Omega)} \cdot I \cdot 
	\bigg(\io |\nabla\varphi|^{2\lambda-2q-2}\bigg)^\frac{1}{2}
  \ee
  with
  \bas
	I:=\Big\| |\nabla\varphi|^{q-1} D^2 \varphi \Big\|_{L^2(\Omega)}.
  \eas
  Likewise, using that $|\Delta\varphi| \le \sqrt{3}|D^2\varphi|$ we can estimate
  \bea{6.4}
	\bigg| - \io \varphi|\nabla\varphi|^{\lambda-2} \Delta\varphi \bigg|
	&\le& \sqrt{3}\|\varphi\|_{L^\infty(\Omega)} \io |\nabla\varphi|^{\lambda-2} |D^2\varphi| \nn\\
	&\le& \sqrt{3} \|\varphi\|_{L^\infty(\Omega)} \cdot I \cdot
	\bigg( \io |\nabla\varphi|^{2\lambda-2q-2} \bigg)^\frac{1}{2}.
  \eea
  Now by the Gagliardo-Nirenberg inequality, there exist $C_1>0$ and $C_2>0$ such that
  \bas
	\bigg( \io |\nabla\varphi|^{2\lambda-2q-2} \bigg)^\frac{1}{2}
	&=& \Big\| |\nabla\varphi|^q \Big\|_{L^\frac{2(\lambda-q-1)}{q}(\Omega)}^\frac{\lambda-q-1}{q} \\
	&\le& C_1 \Big\| \nabla |\nabla \varphi|^q \Big\|_{L^2(\Omega)}^{\frac{\lambda-q-1}{q} \cdot a}
	\Big\| |\nabla\varphi|^q \Big\|_{L^\frac{\lambda}{q}(\Omega)}^{\frac{\lambda-q-1}{q} \cdot (1-a)}
	+ C_1 \Big\| |\nabla\varphi|^q \Big\|_{L^\frac{\lambda}{q}(\Omega)}^\frac{\lambda-q-1}{q} \\
	&\le& C_2 \cdot I^{\frac{\lambda-q-1}{q}\cdot a}
	\cdot \|\nabla\varphi\|_{L^\lambda(\Omega)}^{(\lambda-q-1)(1-a)}
	+ C_1 \|\nabla\varphi\|_{L^\lambda(\Omega)}^{\lambda-q-1}
  \eas
  with
  \bas
	a=\frac{3q(\lambda-2q-2)}{(\lambda-q-1)(6q-\lambda)} \, \in [0,1],
  \eas
  so that combining (\ref{6.2})-(\ref{6.4}) yields $C_3>0$ fulfilling
  \bea{6.5}
	\|\nabla\varphi\|_{L^\lambda(\Omega)}^\lambda
	&\le& C_3 \|\varphi\|_{L^\infty(\Omega)}
	\cdot I^{1+\frac{\lambda-q-1}{q}\cdot a}
	\cdot \|\nabla\varphi\|_{L^\lambda(\Omega)}^{(\lambda-q-1)(1-a)}
	+ C_3 \|\varphi\|_{L^\infty(\Omega)} \cdot I \cdot \|\nabla\varphi\|_{L^\lambda(\Omega)}^{\lambda-q-1} \nn\\
	&=& C_3 \|\varphi\|_{L^\infty(\Omega)}
	\cdot I^\frac{2\lambda-6}{6q-\lambda}
	\cdot \|\nabla\varphi\|_{L^\lambda(\Omega)}^\frac{(4q+1-\lambda)\lambda}{6q-\lambda}
	+ C_3 \|\varphi\|_{L^\infty(\Omega)} \cdot I \cdot \|\nabla\varphi\|_{L^\lambda(\Omega)}^{\lambda-q-1}.
  \eea
  Here we invoke Young's inequality to find $C_4>0$ and $C_5>0$ such that
  \bea{6.6}
	C_3 \|\varphi\|_{L^\infty(\Omega)}
	\cdot I^\frac{2\lambda-6}{6q-\lambda}
	\cdot \|\nabla\varphi\|_{L^\lambda(\Omega)}^\frac{(4q+1-\lambda)\lambda}{6q-\lambda}
	\le \frac{1}{4} \|\nabla\varphi\|_{L^\lambda(\Omega)}^\lambda
	+ C_4 \|\varphi\|_{L^\infty(\Omega)}^\frac{6q-\lambda}{2q-1} \cdot I^\frac{2\lambda-6}{2q-1}
  \eea
  and
  \be{6.7}
	C_3 \|\varphi\|_{L^\infty(\Omega)} \cdot I \cdot \|\nabla\varphi\|_{L^\lambda(\Omega)}^{\lambda-q-1}
	\le \frac{1}{4} \|\nabla\varphi\|_{L^\lambda(\Omega)}^\lambda
	+ C_5 \|\varphi\|_{L^\infty(\Omega)}^\frac{\lambda}{q+1} \cdot I^\frac{\lambda}{q+1},
  \ee
  where by the same token,
  \bea{6.8}
	C_5 \|\varphi\|_{L^\infty(\Omega)}^\frac{\lambda}{q+1} \cdot I^\frac{\lambda}{q+1}
	&=& C_5 \cdot \bigg( \|\varphi\|_{L^\infty(\Omega)}^\frac{6q-\lambda}{2q-1} \cdot I^\frac{2\lambda-6}{2q-1}
	\bigg)^\frac{(2q-1)\lambda}{(q+1)(2\lambda-6)}
	\cdot \|\varphi\|_{L^\infty(\Omega)}^\frac{(3\lambda-6q-6)\lambda}{(q+1)(2\lambda-6)} \nn\\
	&\le& C_5 \|\varphi\|_{L^\infty(\Omega)}^\frac{6q-\lambda}{2q-1} \cdot I^\frac{2\lambda-6}{2q-1}
	+ C_5 \|\varphi\|_{L^\infty(\Omega)}^\lambda.
  \eea
  In consequence, (\ref{6.5})-(\ref{6.8}) prove (\ref{6.1}).
\qed
A first application thereof will appear in the following lemma which estimates the integral
on the right of (\ref{4.1}) under a smallness assumption on the integrability exponent $p$ arrearing therein.
\begin{lem}\label{lem9}
  Let $m\ge 1, q>1$ and $p> \max\{1,m-1\}$ be such that
  \be{9.1}
	p<\Big(2m-\frac{4}{3}\Big) q + m-1.
  \ee
  Then for all $\eta>0$ there exists $C=C(p,q,\eta)>0$ with the property that for all $\eps\in (0,1)$,
  \be{9.4}
	\io \neps^{p+1-m} |\nabla \ceps|^2
	\le \eta \io |\nabla \neps^\frac{p+m-1}{2}|^2 
	+ \eta \io |\nabla \ceps|^{2q-2} |D^2 \ceps|^2
	+ C
	\qquad \mbox{for all } t>0.
  \ee
\end{lem}
\proof
  We apply the H\"older inequality with exponents $\frac{q+1}{q}$ and $q+1$ to obtain
  \bea{9.5}
	\io \neps^{p+1-m} |\nabla \ceps|^2
	&\le& \Big( \io \neps^\frac{(p+1-m)(q+1)}{q} \Big)^\frac{q}{q+1} 
	\cdot \Big( \io |\nabla \ceps|^{2q+2} \Big)^\frac{1}{q+1} \nn\\
	&=& \big\|\neps^\frac{p+m-1}{2}\big\|_{L^\frac{2(p+1-m)(q+1)}{(p+m-1)q}(\Omega)}^\frac{2(p+1-m)}{p+m-1}
	\cdot \|\nabla \ceps\|_{L^{2q+2}(\Omega)}^2
	\qquad \mbox{for all } t>0.
  \eea
  In view of the Gagliardo-Nirenberg inequality (see \cite{win_critexp} for a version involving 
  integrability exponents less than one) and (\ref{mass}), we can find $C_1=C_1(p,q)>0$ and 
  $C_2=C_2(p,q)>0$ such that
  \bea{9.6}
	\big\|\neps^\frac{p+m-1}{2}\big\|_{L^\frac{2(p+1-m)(q+1)}{(p+m-1)q}(\Omega)}^\frac{2(p+1-m)}{p+m-1}
	&\le& C_1 \|\nabla \neps^\frac{p+1-m}{2}\|_{L^2(\Omega)}^{\frac{2(p+1-m)}{p+m-1} \cdot a}
	\|\neps^\frac{p+m-1}{2}\|_{L^\frac{2}{p+m-1}(\Omega)}^{\frac{2(p+1-m)}{p+m-1} \cdot (1-a)} 
	+ C_1 \|\neps^\frac{p+m-1}{2}\|_{L^\frac{2}{p+m-1}(\Omega)}^\frac{2(p+1-m)}{p+m-1} \nn\\
	&\le& C_2 \|\nabla \neps^\frac{p+m-1}{2}\|_{L^2(\Omega)}^{\frac{2(p+1-m)}{p+m-1} \cdot a} + C_2
	\qquad \mbox{for all } t>0,
  \eea
  where $a\in (0,1)$ is determined by the relation
  \be{9.7}
	-\frac{3(p+m-1)q}{2(p+1-m)(q+1)}
	=\Big(1-\frac{3}{2}\Big) \cdot a - \frac{3(p+m-1)}{2} \cdot (1-a).
  \ee
  Here we note that indeed $W^{1,2}(\Omega) \hra L^\frac{2(p+1-m)(q+1)}{(p+m-1)q}(\Omega)$,
  because
  \bas
	6(p+m-1)q - 2(p+1-m)(q+1)
	&=& 6pq + 6(m-1)q - 2pq - 2p + 2(m-1)q + 2(m-1) \\
	&=& (4q-2)p + 8(m-1)q +2(m-1) \\
	&>& 0
  \eas
  and hence $\frac{2(p+1-m)(q+1)}{(p+m-1)q} < 6$.
  On solving (\ref{9.7}) with respect to $a$, we see that (\ref{9.6}) becomes
  \bea{9.8}
	\big\|\neps^\frac{p+m-1}{2}\big\|_{L^\frac{2(p+1-m)(q+1)}{(p+m-1)q}(\Omega)}^\frac{2(p+1-m)}{p+m-1}
	&\le& C_2 \|\nabla \neps^\frac{p+m-1}{2}\|_{L^2(\Omega)}^\frac{6[(p+1-m)(q+1)-q]}{(3p+3m-4)(q+1)}
	+C_2 \nn\\
	&\le& C_3 \cdot \bigg\{ \io |\nabla \neps^\frac{p+m-1}{2}|^2 + 1 \bigg\}^\frac{3[(p+1-m)(q+1)-q]}{(3p+3m-4)(q+1)}
	\quad \mbox{for all } t>0
  \eea
  with some $C_3=C_3(p,q)>0$. \\
  As for the term on the right of (\ref{9.5}) involving $\ceps$, we invoke Lemma \ref{lem6}, 
  which in conjunction with (\ref{cinfty}) provides $C_4=C_4(q)>0$ and $C_5=C_5(q)>0$ satisfying
  \bas
	\|\nabla \ceps\|_{L^{2q+2}(\Omega)}^2
	&\le& C_4 \Big\| |\nabla \ceps|^{q-1} D^2 \ceps\Big\|_{L^2(\Omega)}^\frac{2}{q+1}
	\|\ceps\|_{L^\infty(\Omega)}^\frac{2}{q+1}
	+ C_4 \|\ceps\|_{L^\infty(\Omega)}^2 \\
	&\le& C_5 \cdot \bigg\{ \io |\nabla \ceps|^{2q-2} |D^2 \ceps|^2 + 1 \bigg\}^\frac{1}{q+1}
	\qquad \mbox{for all } t>0.
  \eas
  Together with (\ref{9.8}) and (\ref{9.5}), by Young's inequality this shows that for each $\eta>0$ we can find
  $C_6=C_6(p,q,\eta)>0$ such that
  \bea{9.9}
	\io \neps^{p+1-m} |\nabla \ceps|^2
	&\le& C_3 C_5 \cdot \bigg\| \io |\nabla \neps^\frac{p+m-1}{2}|^2 + 1 \bigg\|^\frac{3[(p+1-m)(q+1)-q]}{(3p+3m-4)(q+1)}
	\cdot \bigg\{ \io |\nabla \ceps|^{2q-2} |D^2 \ceps|^2  + 1 \bigg\}^\frac{1}{q+1} \nn\\
	&\le& \eta \cdot \bigg\{ \io |\nabla \ceps|^{2q-2} |D^2 \ceps|^2 + 1 \bigg\}
	+ C_6 \cdot \bigg\{ \io |\nabla \neps^\frac{p+m-1}{2}|^2 + 1 \bigg\}^\frac{3[(p+1-m)(q+1)-q]}{(3p+3m-4)q}
  \eea
  for all $t>0$.
  Since our hypothesis (\ref{9.1}) warrants that
  \bas
	3[(p+1-m)(q+1)-q] - (3p+3m-4)q
	= 3p-6mq + 4q-3m+3
	<0
  \eas
  and thus
  \bas
	\frac{3[(p+1-m)(q+1)-q]}{(3p+3m-4)q} < 1,
  \eas
  another application of Young's inequality shows that with some $C_7=C_7(p,q,\eta)>0$ we have
  \bas
	C_6 \cdot \bigg\{ \io |\nabla \neps^\frac{p+m-1}{2}|^2 + 1 \bigg\}^\frac{3[(p+1-m)(q+1)-q]}{(3p+3m-4)q}
	\le \eta \cdot \bigg\{ \io |\nabla \neps^\frac{p+m-1}{2}|^2 + 1 \bigg\} + C_7
	\qquad \mbox{for all } t>0.
  \eas
  The claimed inequality (\ref{9.4}) thus results from (\ref{9.9}).
\qed
By pursuing quite a similar strategy, 
under the assumption that $p$ is suitably large as related to $q$ we can estimate the first integral
on the right of (\ref{5.01}), even when enlarged by a further zero-order integral, as follows.
\begin{lem}\label{lem10}
  Let $m\ge 1$ and $q>1$, and suppose that $p>1$ satisfies
  \be{10.1}
	p>\frac{3q-3m+4}{3}.
  \ee
  Then for all $\eta>0$ there exists $C(p,q,\eta)>0$ such that
  \bea{10.4}
	\io \neps^2 |\nabla \ceps|^{2q-2} + \io \neps^2
	\le \eta \io |\nabla \neps^\frac{p+m-1}{2}|^2 + \eta \io |\nabla \ceps|^{2q-2} |D^2 \ceps|^2 
	+ C
	\qquad \mbox{for all } t>0.
  \eea
\end{lem}
\proof
  Once more by the H\"older inequality,
  \bea{10.5}
	\io \neps^2 |\nabla \ceps|^{2q-2}
	&\le& \Big( \io \neps^{q+1} \Big)^\frac{2}{q+1} \cdot 
	\Big(\io |\nabla \ceps|^{2q+2} \Big)^\frac{q-1}{q+1} \nn\\
	&=& \|\neps^\frac{p+m-1}{2}\|_{L^\frac{2(q+1)}{p+m-1}(\Omega)}^\frac{4}{p+m-1} \cdot 
	\|\nabla \ceps\|_{L^{2q+2}(\Omega)}^{2q-2}
	\qquad \mbox{for all } t>0.
  \eea
  Here we observe that (\ref{10.1}) in particular ensures that
  \bas
	6-\frac{2(q+1)}{p+m-1}
	= \frac{6p-2q+6m-8}{p+m-1}
	> \frac{(6q-6m+8)-2q+6m-8}{p+m-1}
	=\frac{4q}{p+m-1}>0
  \eas
  and hence $\frac{2(q+1)}{p+m-1}<6$. 
  We thus may invoke the Gagliardo-Nirenberg inequality, which combined with (\ref{mass}) provides
  $C_1=C_1(p,q)>0$ and $C_2=C_2(p,q)>0$ such that
  \bea{10.6}
	\|\neps^\frac{p+m-1}{2}\|_{L^\frac{2(q+1)}{p+m-1}(\Omega)}^\frac{4}{p+m-1}
	&\le& C_1 \|\nabla \neps^\frac{p+m-1}{2}\|_{L^2(\Omega)}^{\frac{4}{p+m-1} \cdot a}
	\|\neps^\frac{p+m-1}{2}\|_{L^\frac{2}{p+m-1}(\Omega)}^{\frac{4}{p+m-1} \cdot (1-a)}
	+ C_1 \|\neps^\frac{p+m-1}{2}\|_{L^\frac{2}{p+m-1}(\Omega)}^\frac{4}{p+m-1} \nn\\
	&\le& C_2 \cdot \bigg\{ \io |\nabla \neps^\frac{p+m-1}{2}|^2 + 1 \bigg\}^{\frac{2}{p+m-1} \cdot a}
	\qquad \mbox{for all } t>0
  \eea
  with
  \bas
	-\frac{3(p+m-1)}{2(q+1)} = \Big(1-\frac{3}{2}\Big)\cdot a - \frac{3(p+m-1)}{2} \cdot (1-a),
  \eas
  that is, with
  \be{10.7}
	a=\frac{q}{q+1} \cdot \frac{3(p+m-1)}{3p+3m-4} \, \in (0,1).
  \ee
  Likewise, Lemma \ref{lem6} along with (\ref{cinfty}) warrants that
  \bea{10.8}
	\|\nabla \ceps\|_{L^{2q+2}(\Omega)}^{2q-2}
	&\le& C_3 \Big\| |\nabla \ceps|^{q-1} D^2 \ceps \Big\|_{L^2(\Omega)}^\frac{2q-2}{q+1}
	\|\ceps\|_{L^\infty(\Omega)}^\frac{2q-2}{q+1}
	+ C_3 \|\ceps\|_{L^\infty(\Omega)}^{2q-2} \nn\\
	&\le& C_4 \cdot \bigg\{ \io |\nabla \ceps|^{2q-2} |D^2 \ceps|^2 + 1 \bigg\}^\frac{q-1}{q+1}
	\qquad \mbox{for all } t>0
  \eea
  with appropriate positive constants $C_3=C_3(q)$ and $C_4=C_4(q)$.\\
  In light of (\ref{10.6}), (\ref{10.7}) and (\ref{10.8}), we may use Young's inequality in (\ref{10.5}) to see
  that given $\eta>0$ we can find $C_5=C_5(p,q,\eta)>0$ fulfilling
  \bas
	\io \neps^2 |\nabla \ceps|^{2q-2}
	&\le& C_2 C_4 \cdot \bigg\{ \io |\nabla \neps^\frac{p+m-1}{2}|^2 + 1 \bigg\}^\frac{6q}{(q+1)(3p+3m-4)}
	\cdot \bigg\{ \io |\nabla \ceps|^{2q-2} |D^2 \ceps|^2 + 1 \bigg\}^\frac{q-1}{q+1} \\
	&\le& \frac{\eta}{2} \cdot \bigg\{ \io |\nabla \ceps|^{2q-2} |D^2 \ceps|^2 +1 \bigg\}
	+ C_5 \cdot \bigg\{ \io |\nabla \neps^\frac{p+m-1}{2}|^2 + 1 \bigg\}^\frac{3q}{3p+3m-4}
	\quad \mbox{for all } t>0.
  \eas
  Now by (\ref{10.1}) we have
  \bas
	\frac{3q}{3p+3m-4} < \frac{3q}{3\cdot \frac{3q-3m+4}{3} + 3m-4} =1,
  \eas
  so that another application of Young's inequality yields $C_6=C_6(p,q,\eta)>0$ such that
  \bas
	\io \neps^2 |\nabla \ceps|^{2q-2}
	\le \frac{\eta}{2} \io |\nabla \ceps|^{2q-2} |D^2 \ceps|^2 
	+\frac{\eta}{2} \io |\nabla \neps^\frac{p+m-1}{2}|^2
	+ C_6
	\qquad \mbox{for all } t>0.
  \eas
  Since the integral $\io \neps^2$ can be estimated similarly upon a straightforward simplification of the above argument,
  this establishes (\ref{10.4}).
\qed
As for the rightmost integral in (\ref{5.01}), in order to avoid later repetitions in bootstrapping arguments
we find it convenient to estimate this quantity on the basis of a supposedly known bound for $D\ueps$
in $L^\infty((0,\infty);L^r(\Omega))$ for some $r>1$.
In the more favorable case when $r$ is large, under a comparatively mild restriction on $q$ we then
obtain the following.
\begin{lem}\label{lem7}
  Let $m\ge 1$ and $r>\frac{3}{2}$, and suppose that $q\ge r-1$ is such that 
  \be{7.01}
	(4-2r)q \le r-1.
  \ee
  Then for all $\eta>0$ and each $K>0$ there exists $C=C(q,r,\eta,K)>0$ such that if
  for some $\eps\in (0,1)$ and $T>0$ we have
  \be{7.2}
	\|D\ueps(\cdot,t)\|_{L^r(\Omega)} \le K
	\qquad \mbox{for all } t\in (0,T),
  \ee
  then
  \be{7.3}
	\io |\nabla \ceps|^{2q} \cdot |D\ueps|
	\le \eta \io |\nabla \ceps|^{2q-2} |D^2 \ceps|^2 
	+ C
	\qquad \mbox{for all } t\in (0,T).
  \ee
\end{lem}
\proof
  We invoke the H\"older inequality with exponents $\frac{r}{r-1}$ and $r$ to see that
  \bea{7.4}
	\io |\nabla \ceps|^{2q} \cdot |D\ueps|
	&\le& \Big( \io |\nabla \ceps|^\frac{2qr}{r-1} \Big)^\frac{r-1}{r}
	\cdot \Big( \io |D\ueps|^r \Big)^\frac{1}{r} \nn\\
	&\le& K \cdot \|\nabla \ceps\|_{L^\frac{2qr}{r-1}(\Omega)}^{2q}
	\qquad \mbox{for all } t\in (0,T)
  \eea
  due to (\ref{7.2}). 
  Since $q\ge r-1$ ensures that $\lambda:=\frac{2qr}{r-1}$ satisfies
  \bas
	\lambda-(2q+2)
	= \frac{2qr-(2q+2)(r-1)}{r-1} 
	= \frac{2q-2r+2}{r-1}
	= \frac{2(q-r+1)}{r-1} 
	\ge 0,
  \eas
  and since (\ref{7.01}) warrants that
  \bas
	(4q+1)-\lambda
	= \frac{(4q+1)(r-1)-2qr}{r-1}
	= \frac{2qr-4q+r-1}{r-1}
	= \frac{r-1-(4-2r)q}{r-1}
	\ge 0,
  \eas
  we may apply Lemma \ref{lem6} and (\ref{cinfty}) to see that with some $C_1=C_1(q,r)>0$ and $C_2=C_2(q,r)>0$
  we have
  \bea{7.5}
	\|\nabla \ceps\|_{L^\frac{2qr}{r-1}(\Omega)}^{2q}
	&\le& C_1 \Big\| |\nabla \ceps|^{q-1} D^2 \ceps \Big\|_{L^2(\Omega)}^\frac{2q(2\lambda-6)}{(2q-1)\lambda}
	\|\ceps\|_{L^\infty(\Omega)}^\frac{2q(6q-\lambda)}{(2q-1)\lambda}
	+ C_1 \|\ceps\|_{L^\infty(\Omega)}^{2q} \nn\\
	&\le& C_2 \cdot
	\bigg\{ \io |\nabla \ceps|^{2q-2} |D^2 \ceps|^2 \bigg\}^\frac{2q(\lambda-3)}{(2q-1)\lambda} 
	+ C_2
	\qquad \mbox{for all } t\in (0,T),
  \eea
  because clearly $\lambda<6q$. As
  \bas
	\frac{2q(\lambda-3)}{(2q-1)\lambda}
	= \frac{2q(1-\frac{3}{\lambda})}{2q-1}
	= \frac{2q-\frac{3(r-1)}{r}}{2q-1}
	= \frac{2qr-3r+3}{(2q-1)r}
	= \frac{2qr-3r+3}{2qr-r}
	<1
  \eas
  thanks to our assumption $r>\frac{3}{2}$,
  by means of Young's inequality we can easily derive (\ref{7.3}) from (\ref{7.4}) and (\ref{7.5}).
\qed
Also when $r$ is small, however, we can arrive at a similar conclusion, albeit under a slightly stronger
assumption on $q$.
\begin{lem}\label{lem8}
  Let $m\ge 1$, and suppose that $r\in (1,\frac{3}{2}]$ and 
  \bas
	q\in \Big(1,\frac{2r+3}{3}\Big).
  \eas
  Then for all $\eta>0$ and $K>0$ one can find $C=C(q,r,\eta,K)>0$ such that if
  there exist $\eps\in (0,1)$ and $T>0$ fulfilling
  \be{8.2}
	\|D\ueps(\cdot,t)\|_{L^r(\Omega)} \le K
	\qquad \mbox{for all } t\in (0,T),
  \ee
  then
  \be{8.3}
	\io |\nabla \ceps|^{2q} \cdot |D\ueps|
	\le \eta \io |\nabla \ceps|^{2q-2} |D^2 \ceps|^2 
	+ \eta \io |A\ueps|^2
	+ C
	\qquad \mbox{for all } t\in (0,T).
  \ee
\end{lem}
\proof
  By the H\"older inequality,
  \be{8.4}
	\io |\nabla \ceps|^{2q} \cdot |D\ueps|
	\le \|\nabla \ceps\|_{L^{2q+2}(\Omega)}^{2q} \|D\ueps\|_{L^{q+1}(\Omega)}
	\qquad \mbox{for all } t\in (0,T),
  \ee
  where an application of Lemma \ref{lem6} to $\lambda:=2q+2$ in combination with (\ref{cinfty}) 
  yields positive constants $C_1=C_1(q)$ and $C_2=C_2(q)$ satisfying
  \bea{8.5}
	\|\nabla \ceps\|_{L^{2q+2}(\Omega)}^{2q}
	&\le& C_1 \Big\| |\nabla \ceps|^{q-1} D^2 \ceps\Big\|_{L^2(\Omega)}^\frac{2q[2(2q+2)-6]}{(2q-1)(2q+2)}
	\|\ceps\|_{L^\infty(\Omega)}^\frac{2q[6q-(2q+2)]}{(2q-1)(2q+2)}
	+ C_1 \|\ceps\|_{L^\infty(\Omega)}^{2q} \nn\\
	&\le& C_2 \cdot \bigg\{ \io |\nabla \ceps|^{2q-2} |D^2 \ceps|^2 \bigg\}^\frac{q}{q+1}
	+ C_2
	\qquad \mbox{for all } t\in (0,T).
  \eea
  Moreover, using the Gagliardo-Nirenberg inequality and (\ref{8.2}), since $q+1 \le \frac{2r+3}{3}+1<6$ 
  we can find
  $C_3=C_3(q,r)>0$ and $C_4=C_4(q,r)>0$ such that
  \bea{8.6}
	\|D\ueps\|_{L^{q+1}(\Omega)}
	&\le& C_3 \|\ueps\|_{W^{2,2}(\Omega)}^\frac{6(q+1-r)}{(q+1)(6-r)} 
	\|\ueps\|_{W^{1,r}(\Omega)}^\frac{(5-q)r}{(q+1)(6-r)} \nn\\
	&\le& C_4 \|A\ueps\|_{L^2(\Omega)}^\frac{6(q+1-r)}{(q+1)(6-r)} 
	\|D\ueps\|_{L^r(\Omega)}^\frac{(5-q)r}{(q+1)(6-r)} \nn\\
	&\le& C_4 K^\frac{(5-q)r}{(q+1)(6-r)} \|A\ueps\|_{L^2(\Omega)}^\frac{6(q+1-r)}{(q+1)(6-r)}
	\qquad \mbox{for all } t\in (0,T),
  \eea
  because by known regularity estimates for the Stokes operator in bounded domains (see \cite[p.82]{giga_sohr} and the
  references given there)
  and the Poincar\'e inequality we know that $\|A(\cdot)\|_{L^2(\Omega)}$ and $\|D(\cdot)\|_{L^r(\Omega)}$ define
  norms equivalent to $\|\cdot\|_{W^{2,2}(\Omega)}$ and $\|\cdot\|_{W^{1,r}(\Omega)}$, respectively, in $D(A_2)$.\\
  Now given $\eta>0$, we combine (\ref{8.4}), (\ref{8.5}) and (\ref{8.6}) and invoke Young's inequality to find 
  $C_5=C_5(q,r,\eta,K)>0$ such that
  \bea{8.7}
	\io |\nabla \ceps|^{2q} \cdot |D\ueps|
	&\le& C_2 C_4 K^\frac{(5-q)r}{(q+1)(6-r)} \cdot
	\bigg\{ \io |\nabla \ceps|^{2q-2} |D^2 \ceps|^2 \bigg\}^\frac{q}{q+1}
	\cdot \|A\ueps\|_{L^2(\Omega)}^\frac{6(q+1-r)}{(q+1)(6-r)} \nn\\
	& & + C_2 C_4 K^\frac{(5-q)r}{(q+1)(6-r)} \cdot 
	\|A\ueps\|_{L^2(\Omega)}^\frac{6(q+1-r)}{(q+1)(6-r)} \nn\\[1mm]
	&\le& \eta \io |\nabla \ceps|^{2q-2} |D^2 \ceps|^2 
	+ C_5 \|A\ueps\|_{L^2(\Omega)}^\frac{6(q+1-r)}{6-r} \nn\\
	& & + C_2 C_4 K^\frac{(5-q)r}{(q+1)(6-r)} \cdot \|A\ueps\|_{L^2(\Omega)}^\frac{6(q+1-r)}{(q+1)(6-r)}
	\qquad \mbox{for all } t\in (0,T).
  \eea
  Since here our assumption $q<\frac{2r+3}{3}$ entails that
  \bas
	\frac{6(q+1-r)}{(q+1)(6-r)} < \frac{6(q+1-r)}{6-r}<2,
  \eas
  a second application of Young's inequality yields $C_6=C_6(q,r,\eta,K)>0$ fulfilling
  \bas
	C_5 \|A\ueps\|_{L^2(\Omega)}^\frac{6(q+1-r)}{(q+1)(6-r)}
	\le \frac{\eta}{2} \|A\ueps\|_{L^2(\Omega)}^2 + C_6
  \eas
  and
  \bas
	C_2 C_4 K^\frac{(5-q)r}{(q+1)(6-r)} \cdot \|A\ueps\|_{L^2(\Omega)}^\frac{6(q+1-r)}{(q+1)(6-r)}
	\le \frac{\eta}{2} \|A\ueps\|_{L^2(\Omega)}^2 + C_6
  \eas
  for all $t\in (0,T)$. Therefore, (\ref{8.7}) implies (\ref{8.3}).
\qed
\subsection{Combining previous estimates}\label{sect_combine}
Now if $m>\frac{7}{6}$, then the conditions on $p$ in Lemma \ref{lem9} and Lemma \ref{lem10} can be fulfilled
simultaneously for any choice of $q>1$. Thus, resorting to such $m$ allows for combining the above results
to derive an ODI for the functional in (\ref{funct}) which contains a favorable absorptive term.
\begin{lem}\label{lem11}
  Let $m>\frac{7}{6}$. Let $r\ge 1$ and $q>1$ satisfy
  \be{11.1}
	\left\{ \begin{array}{ll}
	q<\frac{2r+3}{3} \qquad & \mbox{if } r\le \frac{3}{2}, \\[1mm]
	(4-2r)q \le r-1 \qquad & \mbox{if } r>\frac{3}{2},
	\end{array} \right.
  \ee
  and assume that $p>\max\{1,m-1\}$ be such that
  \be{11.2}
	\frac{3q-3m+4}{3} < p < \Big(2m-\frac{4}{3}\Big)q + m-1.
  \ee
  Then for all $K>0$ one can find a constant $C=C(p,q,r,K)>0$ such that if for some
  $\eps\in (0,1)$ and $T>0$ we have
  \be{11.4}
	\|D\ueps(\cdot,t)\|_{L^r(\Omega)} \le K
	\qquad \mbox{for all } t\in (0,T),
  \ee
  then
  \bea{11.5}
	\frac{d}{dt} \bigg\{ \io \neps^p + \io |\nabla \ceps|^{2q} + \io |A^\frac{1}{2}\ueps|^2 \bigg\}
	&+& \frac{1}{C} \cdot \bigg\{ \io |\nabla \neps^\frac{p+m-1}{2}|^2 + \io |\nabla \ceps|^{2q-2} |D^2 \ceps|^2 + \io |A\ueps|^2 \bigg\}
	\nn\\[2mm]
	&\le& C
	\qquad \mbox{for all } t\in (0,T).
  \eea
\end{lem}
\proof
  We only need to add suitable multiples of the differential inequalities in Lemma \ref{lem4}, Lemma \ref{lem5} and
  Lemma \ref{lem55} and estimate the terms on the respective right-hand sides by applying Lemma \ref{lem9}, Lemma \ref{lem10}
  and either Lemma \ref{lem8} or Lemma \ref{lem7} with appropriately small $\eta>0$.
\qed
Assuming the above boundedness property of $D\ueps$, upon a further analysis of (\ref{11.5}) we can estimate
$\neps$ in $L^\infty((0,\infty);L^p(\Omega))$ for certain $p\in (1,\infty)$.
\begin{lem}\label{lem12}
  Let $m>\frac{7}{6}$, and assume that $r\ge 1$ and $p>\max\{1,m-1\}$ are such that there exists $q>1$ for which
  (\ref{11.1}) and (\ref{11.2})
  hold. Then for all $K>0$ there exists $C=C(p,q,r,K)>0$ with the property that 
  if $\eps\in (0,1)$ and $T>0$ are such that
  \be{12.2}
	\|D\ueps(\cdot,t)\|_{L^r(\Omega)} \le K
	\qquad \mbox{for all } t\in (0,T),
  \ee
  then we have
  \be{12.31}
	\io \neps^p(\cdot,t) \le C
	\qquad \mbox{for all } t\in (0,T).
  \ee
%
\end{lem}
\proof 
  We only need to derive from (\ref{11.5}) an autonomous ODI for the function
  $y\in C^0([0,T)) \cap C^1((0,T))$ defined by
  \bas
	y(t):=\io \neps^p(\cdot,t) + \io |\nabla \ceps(\cdot,t)|^{2q} + \io |A^\frac{1}{2} \ueps(\cdot,t)|^2,
	\qquad t\in [0,T),
  \eas
  with an appropriate dampening term essentially dominated by
  \bas
	h(t):=\io |\nabla \neps^\frac{p+m-1}{2}(\cdot,t)|^2
	+ \io |\nabla \ceps(\cdot,t)|^{2q-2} |D^2 \ceps(\cdot,t)|^2
	+ \io |A\ueps(\cdot,t)|^2,
	\qquad t\in (0,T).
  \eas
  To this end, we first use Young's inequality, the Poincar\'e inequality and (\ref{mass}) to find positive 
  constants $C_1=C_1(p), C_2=C_2(p)$ and $C_3=C_3(p)$ such that
  \bea{12.5}
	\io \neps^p
	&\le& C_1 \io \neps^{p+m-1} + C_1 \nn\\
	&\le& C_2 \cdot \Big\{ \io |\nabla \neps^\frac{p+m-1}{2}|^2 
	+ \|\neps^\frac{p+m-1}{2}\|_{L^\frac{2}{p+m-1}(\Omega)}^2 \Big\}
	+ C_1 \nn\\
	&\le& C_3 \io |\nabla \neps^\frac{p+m-1}{2}|^2 + C_3
	\qquad \mbox{for all } t\in (0,T).
  \eea
  Next, from Young's inequality, Lemma \ref{lem6} and (\ref{cinfty}) we obtain $C_3=C_3(q)>0, C_4=C_4(q)>0$
  and $C_5=C_5(q)>0$ satisfying
  \bea{12.6}
	\io |\nabla \ceps|^{2q}
	&\le& C_3 \io |\nabla \ceps|^{2q+2} + C_3 \nn\\
	&\le& C_4 \cdot \bigg\{ \Big\| |\nabla \ceps|^{q-1} D^2 \ceps\Big\|_{L^2(\Omega)}^2 	
	\|\ceps\|_{L^\infty(\Omega)}^2
	+ \|\ceps\|_{L^\infty(\Omega)}^{2q+2} \bigg\} + C_3 \nn\\
	&\le& C_5 \io |\nabla \ceps|^{2q-2} |D^2 \ceps|^2 + C_5
	\qquad \mbox{for all } t\in (0,T).
  \eea
  Finally, again by the Poincar\'e inequality we find $C_6>0$ fulfilling
  \be{12.7}
	\io |A^\frac{1}{2}\ueps|^2 = \io |D\ueps|^2 \le C_6 \io |A\ueps|^2
	\qquad \mbox{for all } t\in (0,T),
  \ee
  which when added to (\ref{12.5}) and (\ref{12.6}) shows that
  \bas
	y(t) \le C_7 \cdot h(t) + C_7
	\qquad \mbox{for all } t\in (0,T)
  \eas
  with $C_7:=\max\{C_3,C_5,C_6\}$.\\
  In consequence, an application of Lemma \ref{lem11} yields $C_8=C_8(p,q,r,K)>0$ and 
  $C_9=C_9(p,q,r,K)>0$ such that
  \be{12.8}
	y'(t) + C_8 y(t) +C_8 h(t) \le C_9
	\qquad \mbox{for all } t\in (0,T).
  \ee
  By a comparison argument, this in particular entails that
  \be{12.9}
	y(t) \le C_{10}:=\max \Big\{ y(0) \, , \, \frac{C_9}{C_8} \Big\}
	\qquad \mbox{for all } t\in (0,T),
  \ee
  and thereby proves (\ref{12.31}).	
\qed
Now in light of the mass identity (\ref{mass}), 
a first application of Corollary \ref{cor2} warrants that the hypothesis (\ref{12.2})
in the above lemma is satisfied for some suitably small $r>1$. 
Adjusting the parameter $q$ properly, we thereby arrive at the following result which may be viewed as an
improvement of the regularity property implied by (\ref{mass}), because the number
$5m-\frac{11}{3}$ appearing in (\ref{13.02}) is larger than $1$.
\begin{lem}\label{lem13}
  Let $m>\frac{7}{6}$. Then for all $p>\max\{1,m-1\}$ fulfilling
  \be{13.02}
	p<5m-\frac{11}{3},
  \ee
  one can find $C=C(p)>0$ such that whenever $\eps\in (0,1)$, we have
  \be{13.2}
	\|\neps(\cdot,t)\|_{L^p(\Omega)} \le C
	\qquad \mbox{for all } t>0.
  \ee
\end{lem}
\proof
  We first observe that since $m>1$, we have $m-1<5m-\frac{11}{3}$ and also $\frac{7}{3}-m<5m-\frac{11}{3}$,
  whence without loss of generality we may assume that besides (\ref{13.02}), $p$ satisfies $p>m-1$ and
  \be{13.01}
	p>\frac{7}{3}-m.
  \ee
  Now since $p<5m-\frac{11}{3}$ by (\ref{13.02}), we have $2\cdot(6m-4) > 3p-3m+3$ and hence
  \be{13.3}
	\frac{3(p-m+1)}{6m-4} < 2.
  \ee
  Moreover, our assumption $m>\frac{7}{6}$ ensures that
  \bas
	-18\Big(m-\frac{7}{6}\Big) p < 18 \Big(m-\frac{7}{6}\Big)\Big(m-\frac{1}{3}\Big),
  \eas
  that is,
  \bas
	(21-18m)p < 18 m^2 - 27m + 7,
  \eas
  which is equivalent to the inequality
  \be{13.4}
	\frac{3(p-m+1)}{6m-4} < \frac{3p+3m-4}{3}.
  \ee
  According to (\ref{13.3}), (\ref{13.4}) and the fact that
  \bas
	\frac{3p+3m-4}{3}>1
  \eas
  by (\ref{13.01}), we can now fix $q\in (1,2)$ fulfilling
  \bas
	\frac{3(p-m+1)}{6m-4} < q < \frac{3p+3m-4}{3},
  \eas
  where the left inequality asserts that
  \bas
	p<\Big(2m-\frac{4}{3}\Big)q + m-1,
  \eas
  and the right inequality guarantees that
  \bas
	p>\frac{3q-3m+4}{3},
  \eas
  altogether meaning that (\ref{11.2}) is satisfied.
  Since $q<2$, we can finally pick $r\in [1,\frac{3}{2})$ sufficiently close to $\frac{3}{2}$ such that $r>\frac{3q-3}{2}$,
  so that
  \bas
	q<\frac{2r+3}{3},
  \eas
  ensuring that also (\ref{11.1}) is valid.
  Then in view of (\ref{mass}), Corollary \ref{cor2} asserts that
  \bas
	\|D\ueps(\cdot,t)\|_{L^r(\Omega)} \le C_1
	\qquad \mbox{for all } t>0
  \eas
  with some $C_1>0$, whence according to our choices of $r, q$ and $p$ we may apply Lemma \ref{lem12}
  to find $C_2=C_2(p)>0$ such that
  \bas
	\io \neps^p(\cdot,t) \le C_2
	\qquad \mbox{for all } t>0.
  \eas
  This proves the lemma.
\qed
In a second step, on the basis of the knowledge just gained we may again apply Corollary \ref{cor2}
and once more combine the outcome theoreof with Lemma \ref{lem12} to obtain 
bounds for $\neps$ in any space $L^\infty((0,\infty);L^p(\Omega))$.
\begin{lem}\label{lem14}
  Let $m>\frac{7}{6}$. Then for all $p>1$ 		
  there exists $C=C(p)>0$ such that
  for each $\eps\in (0,1)$ we have
  \be{14.1}
	\|\neps(\cdot,t)\|_{L^p(\Omega)} \le C
	\qquad \mbox{for all } t>0.
  \ee
\end{lem}
\proof
  It is evidently sufficient to show that for any $p_0>\max\{1,m-1\}$ 
  we can find some $p>p_0$ such that (\ref{14.1}) holds with some $C>0$.\\
  For this purpose, given such $p_0$ 	
  we first fix $q>1$ satisfying
  \be{14.5}
	q>\frac{3(p_0+1-m)}{6m-4}
  \ee
  and observe that then since $m>\frac{7}{6}$ we have
  \bas
	3q-3m+4 - \Big[(6m-4)q+3m-3\Big]
	= (7-6m)q + 7-6m < 0
  \eas
  and hence 
  \bas
	\frac{3q-3m+4}{3} < \frac{(6m-4)q+3m-3}{3}.
  \eas
  As (\ref{14.5}) ensures that moreover
  \bas
	\frac{(6m-4)q+3m-3}{3} > \frac{3(p_0+1-m)+3m-3}{3}=p_0,
  \eas
  we can therefore pick some $p>p_0$ fulfilling
  \be{14.6}
	\frac{3q-3m+4}{3} < p < \frac{(6m-4)q+3m-3}{3}.
  \ee
  Now in order to verify (\ref{14.1}) 	
  for these choices of $p$ and $q$, we first use the fact that
  $5m-\frac{11}{3}>\frac{35}{6}-\frac{11}{3}=\frac{13}{6}$ to infer from Lemma \ref{lem13} that there exists $C_1>0$ such that
  \bas
	\|\neps(\cdot,t)\|_{L^\frac{13}{6}(\Omega)} \le C_1
	\qquad \mbox{for all } t>0.
  \eas
  Since $\frac{3\cdot \frac{13}{6}}{3-\frac{13}{6}}=\frac{39}{5}>7$, Corollary \ref{cor2} thereupon yields $C_2>0$ fulfilling
  \bas
	\|D\ueps(\cdot,t)\|_{L^7(\Omega)} \le C_2
	\qquad \mbox{for all } t>0.
  \eas
  As with $r:=7$, the condition $(4-2r)q \le r-1$ in (\ref{11.1}) is trivially satisfied, thanks to (\ref{14.6}) we may thus
  invoke Lemma \ref{lem12} to establish (\ref{14.1}).		
\qed
By means of a Moser-type iteration in conjunction with standard parabolic regularity arguments, 
we can achieve the following boundedness results.
\begin{lem}\label{lem15}
  Let $m>\frac{7}{6}$. Then there exists $C>0$ such that for all $\eps\in (0,1)$ the solution of (\ref{0eps})
  satisfies
  \be{15.1}
	\|\neps(\cdot,t)\|_{L^\infty(\Omega)} \le C
	\qquad \mbox{for all } t>0
  \ee
  and 
  \be{15.2}
	\|\ceps(\cdot,t)\|_{W^{1,\infty}(\Omega)} \le C
	\qquad \mbox{for all } t>0
  \ee
  as well as
  \be{15.3}
	\|\ueps(\cdot,t)\|_{W^{1,\infty}(\Omega)} \le C
	\qquad \mbox{for all } t>0.
  \ee
\end{lem}
\proof
  First, the validity of estimate (\ref{14.1}) in Lemma \ref{lem14} for any $p>3$ allows for an application
  of Corollary \ref{cor2} to $r:=\infty$ to infer that $(D\ueps)_{\eps\in (0,1)}$ 
  is bounded in $L^\infty(\Omega\times (0,\infty))$, and that hence (\ref{15.3}) holds.\\
  Next, using the outcome of Lemma \ref{lem14} with suitably large $p$ and $q$ as a starting point, 
  we may invoke Lemma A.1 in \cite{taowin_JDE} which by means of a Moser-type iteration applied
  to the first equation in (\ref{0eps}) establishes (\ref{15.1}).\\
  Thereupon, (\ref{15.2}) finally can be derived from (\ref{15.1}) and (\ref{15.3}) by well-known
  arguments from parabolic regularity theory for the second equation in (\ref{0eps}) 
  (cf.~the reasoning in \cite[Lemma 4.1]{horstmann_win}, for instance).
\qed
As one further class of a priori estimates,		
let us finally also note straightforward consequences
of Lemma \ref{lem15} for uniform H\"older regularity properties of $\ceps$, $\nabla\ceps$ and $\ueps$.
\begin{lem}\label{lem27}
  Let $m>\frac{7}{6}$.
  Then there exists $\theta\in (0,1)$ such that for some $C>0$ we have
  \be{27.1}
	\|\ceps\|_{C^{\theta,\frac{\theta}{2}}(\bar\Omega \times [t,t+1])} \le C
	\qquad \mbox{for all } t\ge 0,
  \ee
  and such that for each $\tau>0$ we can find $C(\tau)>0$ such that
  \be{27.2}
	\|\nabla \ceps\|_{C^{\theta,\frac{\theta}{2}}(\bar\Omega \times [t,t+1])} \le C
	\qquad \mbox{for all } t\ge \tau.
  \ee
\end{lem}
\proof
  Writing the second equation in the form
  \bas
	c_{\eps t} = \Delta \ceps + g_\eps(x,t),
	\qquad x\in\Omega, \ t>0,
  \eas
  with $g_\eps(x,t):=-\neps f(\ceps) - \ueps \cdot \nabla \ceps$, we immediately obtain both estimates (\ref{27.1})
  and (\ref{27.2}) e.g.~from standard parabolic regularity theory (\cite{LSU}), because $(g_\eps)_{\eps\in (0,1)}$
  is bounded in $L^\infty(\Omega\times (0,\infty))$ according to Lemma \ref{lem15}.
\qed
\begin{lem}\label{lem28}
  Let $m>\frac{7}{6}$.
  Then there exist $\theta\in (0,1)$ and $C>0$ such that
  \be{28.1}
	\|\ueps\|_{C^{\theta,\frac{\theta}{2}}(\bar\Omega \times [t,t+1])} \le C
	\qquad \mbox{for all } t\ge 0.
  \ee
\end{lem}
\proof
  Starting from the variation-of-constants representation
  \be{28.2}
	\ueps(\cdot,t) = e^{-tA} u_0
	+ \int_0^t e^{-(t-s)A} \proj [ \neps(\cdot,s)\nabla \phi] ds,
	\qquad t>0, \ \eps\in (0,1),
  \ee
  we fix $\alpha \in (\frac{3}{4},1)$, so that $D(A_2^\alpha) \hra C^\beta(\bar\Omega)$ for any
  $\beta\in (0,2\alpha-\frac{3}{2})$ (\cite{henry}, \cite{giga1981_the_other}), 
  and apply $A^\alpha$ to both sides of (\ref{28.2}).
  Then performing standard semigroup estimation techniques (\cite{friedman}),
  in view of the boundedness of $(\neps)_{\eps\in (0,1)}$ in $L^\infty(\Omega\times (0,\infty))$ guaranteed by
  Lemma \ref{lem15} we infer the existence of $C_1>0$ and $\beta'>0$ such that for all $\eps\in (0,1)$,
  \bas
	\|A^\alpha \ueps(\cdot,t)\|_{L^2(\Omega)} \le C_1
	\qquad \mbox{for all } t>0
  \eas
  and
  \bas
	\|A^\alpha \ueps(\cdot,t)-A^\alpha \ueps(\cdot,t_0)\|_{L^2(\Omega)} \le C_1 (t-t_0)^{\beta'}
	\qquad \mbox{for all $t_0\ge 0$ and each } t>t_0.
  \eas
  This implies (\ref{28.1}) with some appropriately small $\theta\in (0,1)$.
\qed
\subsection{Some temporally global integrability properties}
Let us next derive three estimates involving integrability over the whole positive time axis. 
They implicitly contain some weak decay properties of the respective integrands, and these properties
will constitute a starting point for our stabilization proof below.\\
The first two of these estimates result from (\ref{0eps}) in a straightforward manner.
\begin{lem}\label{lem21}
  Let $m>\frac{7}{6}$. Then the inequalities
  \be{21.1}
	\int_0^\infty \io \neps f(\ceps) \le \io c_0
  \ee
  and
  \be{21.2}
	\int_0^\infty \io |\nabla\ceps|^2 \le \frac{1}{2} \io c_0^2
  \ee
  are valid for each $\eps\in (0,1)$.
\end{lem}
\proof
  We test the second equation in (\ref{0eps}) by $1$ and $\ceps$ to obtain that for all $\eps\in (0,1)$ and $t>0$,
  \bas
	\io \ceps(\cdot,t) + \int_0^t \io \neps f(\ceps) = \io c_0
  \eas
  and
  \bas
	\frac{1}{2} \io \ceps^2(\cdot,t) + \int_0^t \io |\nabla\ceps|^2 + \int_0^t \io \neps \ceps f(\ceps)
	= \frac{1}{2} \io c_0^2,
  \eas
  respectively. From these identities, (\ref{21.1}) and (\ref{21.2}) immediately follow.
\qed
A corresponding spatio-temporal integrability property of $\nabla\neps$ can be obtained from (\ref{4.1})
upon using (\ref{21.2}) along with the $L^\infty$ bound for $\neps$ from Lemma \ref{lem15}.
\begin{lem}\label{lem22}
  Assume that $m>\frac{7}{6}$, and that $p>1$ is such that $p\ge m-1$. Then there exists $C>0$ with the property that
  \be{22.1}
	\int_0^\infty \io |\nabla \neps^\frac{p+m-1}{2}|^2 \le C
  \ee
  for all $\eps\in (0,1)$.
\end{lem}
\proof
  Due to Lemma \ref{lem15}, there exists $C_1>0$ such that for all $\eps\in (0,1)$ we have
  $\neps \le C_1$ in $\Omega\times (0,\infty)$.
  Since $p\ge m-1$, we can thereby estimate the integral on the right of (\ref{4.1}) according to
  \bas
	\io \neps^{p+1-m}|\nabla\ceps|^2
	\le C_1^{p+1-m} \io |\nabla\ceps|^2
	\qquad \mbox{for all $t>0$ and } \eps\in (0,1).
  \eas
  Hence, an integration of (\ref{4.1}) shows that
  \bas
	\frac{2(p-1)\kd}{(p+m-1)^2} \int_0^t \io |\nabla \neps^\frac{p+m-1}{2}|^2
	\le \frac{1}{p} \io n_0^p
	+ \frac{(p-1) S_1^2}{2\kd} \cdot C_1^{p+1-m} \int_0^t \io |\nabla \ceps|^2
  \eas
  for any such $t$ and $\eps$, from which (\ref{22.1}) results by an application of Lemma \ref{lem21}.
\qed
\subsection{Regularity properties of time derivatives}
In order to pass to the limit in the first equation in (\ref{0eps}), we shall need an appropriate boundedness
property of the time derivatives of certain powers of $\neps$. 
On time intervals of a fixed finite length, this can be achieved in a straightforward way by making use of
the a priori bounds derived so far.
\begin{lem}\label{lem23}
  Suppose that $m>\frac{7}{6}$, and let $\gamma>m$ satisfy $\gamma \ge 2(m-1)$.
  Then for all $T>0$ there exists $C(T)>0$ such that
  \be{23.1}
	\int_0^T \|\partial_t \neps^\gamma(\cdot,t)\|_{(W_0^{3,2}(\Omega))^\star)} dt \le C(T)
	\qquad \mbox{for all } \eps\in (0,1).
  \ee
\end{lem}
\proof
  On differentiation and integration by parts in (\ref{0eps}), we see that for each fixed $\psi\in C_0^\infty(\Omega)$
  we have
  \bea{23.2}
	\frac{1}{\gamma} \io \partial_t \neps^\gamma (\cdot,t)\cdot\psi
	&=& \io \neps^{\gamma-1} \cdot \bigg\{ \nabla \cdot \Big(D_\eps(\neps)\nabla \neps\Big) 
	- \nabla \cdot \Big(\neps S_\eps(x,\neps,\ceps)\cdot\nabla\ceps\Big) - \ueps \cdot \nabla \neps \bigg\} 
	\cdot \psi \nn\\[1mm]
	&=& - (\gamma-1) \io \neps^{\gamma-2} D_\eps(\neps) |\nabla \neps|^2 \psi
	- \io \neps^{\gamma-1} D_\eps(\neps) \nabla \neps \cdot \nabla \psi \nn\\
	& & + (\gamma-1) \io \neps^{\gamma-1} \nabla\neps \cdot \Big(S_\eps(x,\neps,\ceps) \cdot\nabla \ceps\Big) \psi
	+ \io \neps^\gamma \Big(S_\eps(x,\neps,\ceps)\cdot\nabla\ceps\Big) \cdot\nabla\psi \nn\\
	& & + \frac{1}{\gamma} \io \neps^\gamma \ueps \cdot \nabla \psi
	\qquad \mbox{for all } t>0.
  \eea
  In order to estimate the integrals on the right appropriately, we first apply Lemma \ref{lem15} to fix positive
  constants $C_1, C_2$ and $C_3$ such that
  \be{23.3}
	\neps \le C_1, \quad |\nabla \ceps| \le C_2
	\quad \mbox{and} \quad
	|\ueps| \le C_3
	\qquad \mbox{in } \Omega\times (0,\infty)
	\qquad \mbox{for all } \eps\in (0,1),
  \ee
  whence due to the fact that $D_\eps \le D+2\eps$ in $(0,\infty)$, we also have
  \be{23.4}
	D_\eps(\neps) \le C_4 := \|D_0\|_{L^\infty((0,C_1))} + 2
	\qquad \mbox{in } \Omega\times (0,\infty)
	\qquad \mbox{for all } \eps\in (0,1).
  \ee
  Moreover, since $\gamma>m$ and $\gamma\ge 2(m-1)$, the number $p:=\gamma-m+1$ satisfies $p>1$ and $p\ge m-1$,
  so that Lemma \ref{lem22} becomes applicable so as to yield $C_5>0$ fulfilling
  \be{23.5}
	\int_0^\infty \io \neps^{\gamma-2} |\nabla \neps|^2 
	= \int_0^\infty \io \neps^{p+m-3} |\nabla\neps|^2
	\le C_5
	\qquad \mbox{for all } \eps\in (0,1).
  \ee
  Now using (\ref{23.3}), (\ref{23.4}) and Young's inequality, in (\ref{23.2}) we find that
  \be{23.6}
	\bigg| -(\gamma-1) \io \neps^{\gamma-2} D_\eps(\neps) |\nabla\neps|^2 \psi \bigg|
	\le (\gamma-1) C_4 \cdot \Big( \io \neps^{\gamma-2} |\nabla \neps|^2 \bigg) \cdot \|\psi\|_{L^\infty(\Omega)}
  \ee
  as well as
  \bea{23.7}
	\bigg|-\io \neps^{\gamma-1} D_\eps(\neps) \nabla\neps \cdot \nabla \psi \bigg|
	&\le& C_4 \cdot \bigg( \io \neps^{\gamma-1} |\nabla \neps| \bigg) \cdot \|\nabla\psi\|_{L^\infty(\Omega)} \nn\\
	&\le& C_4 \cdot \bigg\{ \io \neps^{\gamma-2} |\nabla\neps|^2 + \io \neps^\gamma \bigg\} 
	\cdot \|\nabla\psi\|_{L^\infty(\Omega)} \nn\\
	&\le& \bigg\{ C_4 \io \neps^{\gamma-2} |\nabla\neps|^2 + C_4 C_1^\gamma |\Omega| \bigg\} 
	\cdot \|\nabla\psi\|_{L^\infty(\Omega)}
  \eea
  and, similarly,
  \bea{23.8}
	& & \hspace*{-40mm}
	\bigg| (\gamma-1) \io \neps^{\gamma-1} \nabla\neps \cdot \Big(S_\eps(x,\neps,\ceps)\cdot\nabla \ceps \Big) \psi \bigg|
	\nn\\
	&\le& (\gamma-1) \cdot \bigg( \io \neps^{\gamma-1}|\nabla\neps| \bigg) \cdot S_1 C_2 \|\psi\|_{L^\infty(\Omega)} \nn\\
	&\le& (\gamma-1) S_1 C_2 \cdot \bigg\{ \io \neps^{\gamma-2} |\nabla\neps|^2 + C_1^\gamma |\Omega| \bigg\}
	\cdot \|\psi\|_{L^\infty(\Omega)}
  \eea 
  whereas by means of (\ref{23.3}) and (\ref{cinfty}) we can estimate
  \be{23.9}
	\bigg| \io \neps^\gamma \Big(S_\eps(x,\neps,\ceps)\cdot\nabla\ceps\Big)\cdot\nabla\psi \bigg|
	\le C_1^\gamma S_1 C_2 |\Omega| \|\nabla\psi\|_{L^\infty(\Omega)}
  \ee
  and
  \be{23.10}
	\bigg| \frac{1}{\gamma} \io \neps^\gamma \ueps \cdot \nabla \psi \bigg|
	\le \frac{1}{\gamma} C_1^\gamma C_3 |\Omega| \|\nabla\psi\|_{L^\infty(\Omega)}
  \ee
  for all $\eps\in (0,1)$.\\
  As in the considered three-dimensional setting we have $W_0^{3,2}(\Omega) \hra W^{1,\infty}(\Omega)$, collecting
  (\ref{23.6})-(\ref{23.10}) we infer the existence of $C_6>0$ such that
  \bas
	\|\partial_t \neps^\gamma(\cdot,t)\|_{(W_0^{3,2}(\Omega))^\star}
	\le C_6 \cdot \bigg\{ \io \neps^{\gamma-2} |\nabla\neps|^2 + 1 \bigg\}
	\qquad \mbox{for all $t>0$ and any } \eps \in (0,1).
  \eas
  According to (\ref{23.5}), for each $T>0$ we therefore have
  \bas
	\int_0^T \|\partial_t \neps^\gamma(\cdot,t)\|_{(W_0^{3,2}(\Omega))^\star} dt \le
	C_5 C_6 + C_6 T
	\qquad \mbox{for all } \eps\in (0,1),
  \eas
  which proves (\ref{23.1}).
\qed
In proving that the limit function $n$ gained below approaches its spatial mean not only along certain sequences
of times but in fact along the entire net $t\to\infty$, we will rely on an additional regularity estimate
for $n_{\eps t}$ which, in contrast to that in Lemma \ref{lem23}, is uniform with respect to time.
\begin{lem}\label{lem24}
  Let $m>\frac{7}{6}$. Then there exists $C>0$ such that 
  \be{24.1}
	\|n_{\eps t}(\cdot,t)\|_{(W_0^{2,2}(\Omega))^\star} \le C
	\qquad \mbox{for all $t>0$ and } \eps\in (0,1).
  \ee
  In particular,
  \be{24.11}
	\|\neps(\cdot,t)-\neps(\cdot,s)\|_{(W_0^{2,2}(\Omega))^\star} \le C |t-s|
	\qquad \mbox{for all $t\ge 0, s\ge 0$ and } \eps\in (0,1).
  \ee
\end{lem}
\proof
  We fix $\psi\in C_0^\infty(\Omega)$ and multiply the first equation in (\ref{0eps}) by $\psi$. Integrating by parts
  we find that
  \bea{24.2}
	\io n_{\eps t}(\cdot,t)\psi
	&=& \io \nabla \cdot (D_\eps(\neps)\nabla\neps)\psi
	- \io \nabla \cdot \Big( \neps S_\eps(x,\neps,\ceps)\cdot\nabla\ceps\Big)\psi
	- \io (\ueps\cdot \nabla\neps) \psi \nn\\
	&=& \io H_\eps(\neps) \Delta \psi
	+ \io \neps \Big(S_\eps(x,\neps,\ceps)\cdot\nabla\ceps \Big)\cdot\nabla\psi
	+ \io \neps \ueps \cdot\nabla \psi
  \eea
  for all $t>0$, where we have set $H_\eps(s):=\int_0^s D_\eps(\sigma)d\sigma$ for $s\ge 0$.
  Here since by Lemma \ref{lem15} we can find $C_1>0$ such that
  $\neps \le C_1$ in $\Omega\times (0,\infty)$ for all $\eps\in (0,1)$,
  recalling that $D_\eps \le D+2\eps$ on $(0,\infty)$ we can estimate
  \bas
	H_\eps(\neps) \le C_2:=C_1 \cdot \Big(\|D\|_{L^\infty((0,C_1))} + 2 \Big)
	\qquad \mbox{in $\Omega\times (0,\infty)$ \quad for all } \eps\in (0,1).
  \eas
  If moreover we invoke Lemma \ref{lem15} once again to pick $C_3>0$ and $C_4>0$ fulfilling
  \bas
	|\nabla \ceps| \le C_3
	\quad \mbox{and} \quad
	|\ueps| \le C_4
	\qquad \mbox{in $\Omega\times (0,\infty)$ \quad for all } \eps\in (0,1),
  \eas
  then from (\ref{24.2}), (\ref{S2}) and (\ref{cinfty}) we can derive the inequality
  \bas
	\bigg| \io n_{\eps t}(\cdot,t) \cdot\psi \bigg|
	\le C_2 \io |\Delta\psi|
	+ C_1 C_3 S_1 \io |\nabla \psi|
	+ C_1 C_4 \io |\nabla\psi|
	\qquad \mbox{for all $t>0$ and } \eps\in (0,1),
  \eas
  where $S_1:=S_0(\|c_0\|_{L^\infty(\Omega)}$.
  This readily establishes (\ref{24.1}) and thus also (\ref{24.11}).
\qed
\mysection{Passing to the limit. Proof of Theorem \ref{theo16}}
Our generalized solution concept reads as follows.
\begin{defi}\label{defi_weak}
  Suppose that $(n_0,c_0,u_0)$ satisfies (\ref{init}), and let $T>0$.
  Then by a {\em weak solution} of (\ref{0}) in $\Omega\times (0,T)$ we mean a triple of functions
  \bea{reg_w1}
	& & n \in L^1_{loc}(\bar\Omega \times [0,T)), \nn\\
	& & c \in L^\infty_{loc}(\bar\Omega \times [0,T)) \cap L^1_{loc}([0,T);W^{1,1}(\Omega)), \nn\\
	& & u \in L^1_{loc}([0,T);W^{1,1}(\Omega)),
  \eea
  such that $n\ge 0$ and $c\ge 0$ in $\Omega\times (0,T)$ and
  \bea{reg_w2}
	H(n), \ n|\nabla c| \ \mbox{and} \ n|u|
	\ \mbox{ belong to } L^1_{loc}(\bar\Omega\times [0,T)),
  \eea
  that $\nabla \cdot u=0$ in the distributional sense in $\Omega\times (0,T)$, and such that
  \be{w1}
	- \int_0^T \io n\varphi_t - \io n_0 \varphi(\cdot,0)
	= \int_0^T \io H(n) \Delta \varphi
	+ \int_0^T \io n \Big(S(x,n,c)\cdot\nabla c\Big) \cdot \nabla \varphi
	+ \int_0^T \io nu\cdot \nabla \varphi
  \ee
  for all $\varphi\in C_0^\infty(\bar\Omega\times [0,T))$ fulfilling $\frac{\partial\varphi}{\partial\nu}=0$ on
  $\partial\Omega\times (0,T)$, that
  \be{w2}
	- \int_0^T \io c\varphi_t - \io c_0 \varphi(\cdot,0)
	= - \int_0^T \io \nabla c \cdot \nabla \varphi
	- \int_0^T \io nf(c)\varphi
	+ \int_0^T \io cu\cdot \nabla \varphi
  \ee
  for all $\varphi\in C_0^\infty(\bar\Omega\times [0,T))$, and that
  \be{w3}
	- \int_0^T \io u\cdot \varphi_t - \io u_0\cdot \varphi(\cdot,0)
	= - \int_0^T \io \nabla u \cdot \nabla \varphi
	+ \int_0^T \io n\nabla \phi \cdot \varphi
  \ee
  for all $\varphi \in C_0^\infty(\Omega\times [0,T);\R^3)$ such that $\nabla \varphi\equiv 0$ in $\Omega\times (0,T)$,
  where we have set
  \bas
	H(s):=\int_0^s D(\sigma)d\sigma
	\qquad \mbox{for } s\ge 0.
  \eas
  If $(n,c,u) : \Omega\times (0,\infty)\to \R^5$ is a weak solution of (\ref{0}) in $\Omega\times (0,T)$ for all $T>0$,
  then we call $(n,c,u)$ a {\em global weak solution} of (\ref{0}).
\end{defi}
In this framework, (\ref{0}) is indeed globally solvable. This can be seen by making use of the above a priori
estimates and extracting suitable subsequences in a standard manner.
\begin{lem}\label{lem26}
  Let $m>\frac{7}{6}$. Then there exists $(\eps_j)_{j\in\N} \subset (0,1)$ such that $\eps_j\searrow 0$ as $j\to\infty$ and
  that
  \begin{eqnarray}
	& & \neps \to n
	\qquad \mbox{a.e.~in } \Omega\times (0,\infty),
	\qquad \label{26.1} \\
	& & \neps \wsto n
	\qquad \mbox{in } L^\infty(\Omega\times (0,\infty)),
	\qquad \label{26.2} \\
	& & \neps \to n
	\qquad \mbox{in } C^0_{loc}\Big([0,\infty); (W_0^{2,2}(\Omega))^\star\Big),
	\qquad \label{26.3} \\
	& & \ceps \to c
	\qquad \mbox{in } C^0_{loc}(\bar\Omega\times [0,\infty)),
	\qquad \label{26.4} \\
	& & \nabla \ceps \to \nabla c
	\qquad \mbox{in } C^0_{loc}(\bar\Omega\times (0,\infty)),
	\qquad \label{26.5} \\
	& & \nabla \ceps \wsto \nabla c
	\qquad \mbox{in } L^\infty(\Omega\times (0,\infty)),
	\qquad \label{26.6} \\
	& & \ueps \to u
	\qquad  \mbox{in } C^0_{loc}(\bar\Omega\times [0,\infty)) \qquad \mbox{and}
	\qquad \label{26.7} \\
	& & D \ueps \wsto Du
	\qquad \mbox{in } L^\infty(\Omega\times (0,\infty))
	\qquad \label{26.8}
  \eea
  with some triple $(n,c,u)$ which is a global weak solution of (\ref{0}) in the sense of Definition \ref{defi_weak}.
  Moreover, $n$ satisfies
  \be{26.9}
	n\in C^0_{w-\star}([0,\infty);L^\infty(\Omega))
  \ee
  as well as
  \be{mass_n}
	\io n(x,t)dx = \io n_0(x)dx
	\qquad \mbox{for all } t>0.
  \ee
\end{lem}
\proof
  In view of Lemma \ref{lem27} and Lemma \ref{lem28}, the Arzel\`a-Ascoli theorem along with a standard extraction procedure
  yields a sequence $(\eps_j)_{j\in\N}\subset (0,1)$ with $\eps_j \searrow 0$ as $j\to\infty$ such that
  (\ref{26.4}), (\ref{26.5}) and (\ref{26.7}) hold with some limit functions $c$ and $u$ belonging to the indicated spaces.
  Passing to a subsequence if necessary, by means of Lemma \ref{lem15} we can achieve that for some 
  $n\in L^\infty(\Omega\times (0,\infty))$ we moreover have (\ref{26.2}), (\ref{26.6}) and (\ref{26.8}).\\
  We next fix $\gamma>m$ such that $\gamma \ge 2(m-1)$ and combine Lemma \ref{lem23} with the estimate asserted by
  Lemma \ref{lem22} for $p:=2\gamma-m+1$ to see that for each $T>0$, $(\eps^\gamma)_{\eps\in (0,1)}$ is bounded
  in $L^2((0,T);W^{1,2}(\Omega))$ with $(\partial_t \neps^\gamma)_{\eps\in (0,1)}$ being bounded in
  $L^1((0,T);(W_0^{3,2}(\Omega))^\star)$. Therefore, an Aubin-Lions lemma (\cite{temam}) applies to yield strong
  precompactness of $(\neps^\gamma)_{\eps\in (0,1)}$ in $L^2(\Omega\times (0,T))$, whence along a suitable subsequence
  we have $\neps^\gamma \to z^\gamma$ and hence $\neps\to z$ a.e.~in $\Omega\times (0,\infty)$ for some nonnegative measurable
  $z:\Omega\times (0,\infty)\to \R$. By Egorov's theorem, we know that necessarily $z=n$, so that (\ref{26.1}) follows.\\
  Finally, as the embedding $L^\infty(\Omega) \hra (W_0^{2,2}(\Omega))^\star$ is compact, the Arzel\`a-Ascoli once more
  applies to say that the equicontinuity property (\ref{24.11}) together with the boundedness of $(\neps)_{\eps\in (0,1)}$
  in $C^0([0,\infty);L^\infty(\Omega))$ ensures that (\ref{26.3}) holds after a further extraction of an adequate
  subsequence.\abs
  The additional regularity property (\ref{26.9}) thereafter is a consequence of (\ref{26.3}) and the fact that
  $C_1:=\|n\|_{L^\infty(\Omega\times (0,\infty))}$ is finite: First, from the latter two properties it can readily be seen
  that $n(\cdot,t)\in L^\infty(\Omega)$ actually for all $t\in [0,\infty)$, with $\|n(\cdot,t)\|_{L^\infty(\Omega)} \le C_1$
  for all $t\ge 0$. Thus, given any $t_0\ge 0$ and $(t_j)_{j\in\N}\subset [0,\infty)$ such that $t_j\to t_0$ as $j\to\infty$
  we know that $(n(\cdot,t_j))_{j\in\N}$ is bounded in $L^\infty(\Omega)$, and that for all $\psi\in C_0^\infty(\Omega)$
  we have
  $\io n(\cdot,t_j)\psi \to \io n(\cdot,t_0)\psi$ as $j\to\infty$ by (\ref{26.3}). 
  By density of $C_0^\infty(\Omega)$ in $L^1(\Omega)$, this proves that indeed $n(\cdot,t_j) \wsto n(\cdot,t_0)$ as 
  $j\to\infty$.\abs
  Now the verification of the claimed weak solution property of $(n,c,u)$ is straightforward:
  Whereas the nonnegativity of $n$ and $c$ and the integrability requirements in (\ref{reg_w1}) and (\ref{reg_w2})
  are immediate from (\ref{26.1}), (\ref{26.2}), (\ref{26.4}), (\ref{26.6}) and (\ref{26.7}), the integral identities
  (\ref{w1}), (\ref{w2}) and (\ref{w3}) can be derived by standard arguments from the corresponding weak formulations in 
  the approximate system (\ref{0eps}) upon letting $\eps=\eps_j\searrow 0$ and using (\ref{26.1}) and (\ref{26.2}) as 
  well as (\ref{26.4})-(\ref{26.8}).
\qed
\proofc of Theorem \ref{theo16}.\quad
  We only need to combine Lemma \ref{lem26} with Lemma \ref{lem15}.
\qed
\mysection{Large time behavior. Proof of Theorem \ref{theo166}}\label{sect_conv}
In this section we shall assume that $D, S$ and $f$ satisfy the assumptions in Theorem \ref{theo166},
and we will establish the convergence statements therein separately for the solution components $n,c$ and $u$.\abs
Here proving stabilization of $n$ will require a comparatively subtle reasoning, which is due to the fact
that our knowledge on compactness properties of $(n(\cdot,t))_{t>0}$, and of temporal continuity features of $n$,
is rather limited.
The core of the following argument lies in an appropriate combination of the decay information implied
by Lemma \ref{lem22} with the continuity property contained in Lemma \ref{lem24}.
The nonlinearity of diffusion in the first equation in (\ref{0}), reflected in the appearance of a nontrivial function 
of $n$ in the integral in (\ref{22.1}), requires the use of certain powers of $n$ in the following proof. 
\begin{lem}\label{lem30}
  Let $m>\frac{7}{6}$. Then with $(n,c,u)$ as given by Theorem \ref{theo16}, we have
  \be{30.1}
	n(\cdot,t)\wsto \onz
	\quad \mbox{in } L^\infty(\Omega)
	\qquad \mbox{as } t\to\infty.
  \ee
\end{lem}
\proof
  Let us assume that the conclusion of the lemma does not hold. Then we can find a sequence 
  $(t_j)_{j\in\N} \subset (0,\infty)$ such that $t_j\to\infty$ as $j\to\infty$, and such that for some 
  $\widetilde\psi\in L^1(\Omega)$ we have
  \be{30.2}
	\io n(x,t_j)\widetilde \psi(x)dx
	-\io \onz \widetilde \psi(x) dx \ge C_1
	\qquad \mbox{for all } j\in\N
  \ee
  with some $C_1>0$. To exploit his appropriately,
  according to Lemma \ref{lem15} we take $C_2>0$ fulfilling
  \be{30.22}
	n(x,t)\le C_2
	\qquad \mbox{for a.e.~$(x,t)\in \Omega\times (0,\infty$),}
  \ee
  and then use the density of $C_0^\infty(\Omega)$ in $L^1(\Omega)$ in choosing $\psi\in C_0^\infty(\Omega)$ such that
  $\|\psi-\widetilde\psi\|_{L^1(\Omega)} \le \frac{C_1}{4C_2}$, so that by (\ref{30.2}),
  \bea{30.3}
	\io n(x,t_j)\psi(x)dx - \io \onz \psi(x)dx
	&\ge& \io n(x,t_j)\widetilde\psi(x)dx - \io \onz \widetilde{\psi}(x)dx \nn\\
	& & - \Big\{ \|n(\cdot,t_j)\|_{L^\infty(\Omega)} + \onz \Big\} \cdot \|\psi-\widetilde\psi\|_{L^1(\Omega)} \nn\\[2mm]
	&\ge& \frac{C_1}{2}
	\qquad \mbox{for all } j\in\N.
  \eea
  Now since by Lemma \ref{lem24} there exists $C_3>0$ such that for all $\eps\in (0,1)$ we have
  \bas
	\|\neps(\cdot,t)-\neps(\cdot,s)\|_{(W_0^{2,2}(\Omega))^\star} \le C_3 |t-s|
	\qquad \mbox{for all $t\ge 0$ and $s\ge 0$}, 
  \eas
  recalling the convergence property (\ref{26.3}) from Lemma \ref{lem26} we see that
  \bas
	\|n(\cdot,t)-n(\cdot,s)\|_{(W_0^{2,2}(\Omega))^\star} \le C_3 |t-s|
	\qquad \mbox{for all $t\ge 0$ and $s\ge 0$}, 
  \eas
  In particular, this implies that if we let $\tau\in (0,1)$ be such that
  \bas
	\tau \le \frac{C_1}{4C_3 \|\psi\|_{W_0^{2,2}(\Omega)}},
  \eas
  then for all $j\in\N$ and each $t\in (t_j,t_j+\tau)$ we have
  \bas
	\bigg| \io n(x,t_j) \psi(x)dx - \io n(x,t) \psi(x)dx \bigg|
	&\le& \|n(\cdot,t_j)-n(\cdot,t)\|_{(W_0^{2,2}(\Omega))^\star} \cdot \|\psi\|_{W_0^{2,2}(\Omega)} \\
	&\le& C_3 |t_j-t| \cdot \|\psi\|_{W_0^{2,2}(\Omega)} \\[1mm]
	&\le& \frac{C_1}{4},
  \eas
  and hence (\ref{30.3}) entails that
  \be{30.4}
	\io n(x,t)\psi(x)dx - \io \onz \psi(x) dx
	\ge \frac{C_1}{4}
	\qquad \mbox{for all $t\in (t_j,t_j+\tau)$ and each } j\in\N.
  \ee
  To see that this contradicts the outcome of Lemma \ref{lem22}, we fix any $p>1$ such that $p\ge m-1$ and $p>3-m$,
  and abbreviate $\gamma:=\frac{p+m-1}{2}>1$. Then taking a Poincar\'e constant $C_4>0$ such that
  \bas
	\io |\varphi(x) -\overline{\varphi}|^2 dx \le C_4 \io |\nabla\varphi|^2
	\qquad \mbox{for all } \varphi\in W^{1,2}(\Omega),
  \eas
  again with $\overline{\varphi}:=\frac{1}{|\Omega|}\io \varphi$, from Lemma \ref{lem22} we obtain $C_5>0$ such that
  \bea{30.5}
	\int_0^\infty \io |\neps^\gamma(x,t)-a_\eps^\gamma(t)|^2 dxdt
	&\le& C_4 \int_0^\infty \io |\nabla \neps^\gamma|^2 \nn\\[2mm]
	&\le& C_5
	\qquad \mbox{for all } \eps\in (0,1),
  \eea
  where
  \bas
	a_\eps(t):=\Big( \overline{\neps^\gamma(\cdot,t)} \Big)^\frac{1}{\gamma}
	= \bigg\{ \frac{1}{|\Omega|} \io \neps^\gamma(x,t)dx \bigg\}^\frac{1}{\gamma}
	\qquad \mbox{for $\eps\in (0,1)$ and } t>0.
  \eas
  Here since from Lemma \ref{lem26} and the Tonelli theorem we know that for a.e.~$t>0$ we have $\neps(\cdot,t)\to 
  n(\cdot,t)$ a.e.~in $\Omega$ as $\eps=\eps_j\searrow 0$, in view of the uniform boundedness of $(\neps)_{\eps\in (0,1)}$
  in $L^\infty(\Omega\times (0,\infty))$ we may apply the dominated convergence theorem to infer that
  \be{30.6}
	a_\eps(t)\to a(t)
	\qquad \mbox{for a.e.~$t>0$}
  \ee
  as $\eps=\eps_j\searrow 0$, where 
  \be{30.66}
	a(t):=\Big( \overline{n^\gamma(\cdot,t)} \Big)^\frac{1}{\gamma}
	= \bigg\{ \frac{1}{|\Omega|} \io n^\gamma(x,t)dx \bigg\}^\frac{1}{\gamma}
	\qquad \mbox{for } t>0.
  \ee
  Again using that $\neps\to n$ a.e.~in $\Omega\times (0,\infty)$ as $\eps=\eps_j\searrow 0$, by means of Fatou's lemma
  we can derive from (\ref{30.5}) and (\ref{30.6}) that
  \be{30.7}
	\int_0^\infty \io |n^\gamma(x,t)-a^\gamma(t)|^2 dxdt \le C_5.
  \ee
  Thanks to the fact that $\gamma>1$ ensures the validity of the elementary inequality 
  \bas
	\frac{\xi^\gamma-\eta^\gamma}{\xi-\eta} \ge \eta^{\gamma-1}
	\qquad \mbox{for all $\xi\ge 0$ and $\eta\ge0$ such that $\xi\ne\eta$,}
  \eas
  and since by the H\"older inequality, (\ref{mass_n}) and (\ref{30.66}) we have
  \bas
	\onz
	= \frac{1}{|\Omega|} \io n(x,t)dx
	\le \frac{1}{|\Omega|} \cdot \bigg(\io n^\gamma(x,t)dx \bigg)^\frac{1}{\gamma} \cdot |\Omega|^\frac{\gamma-1}{\gamma}
	= a(t)
	\qquad \mbox{for a.e.~$t>0$,}
  \eas
  on the left of (\ref{30.7}) we can estimate
  \bas
	\io |n^\gamma(x,t)-a^\gamma(t)|^2 dx
	&\ge& a^{2\gamma-2}(t) \cdot \io |n(x,t)-a(t)|^2 dx \\
	&\ge& \onz^{2\gamma-2} \cdot  \io |n(x,t)-a(t)|^2 dx
	\qquad \mbox{for a.e.~$t>0$.}
  \eas
  Therefore, (\ref{30.7}) yields
  \be{30.8}
	\int_0^\infty \io |n(x,t)-a(t)|^2 dxdt \le C_6:=\frac{C_5}{\onz^{2\gamma-2}}.
  \ee
  We now introduce
  \bas
	n_j(x,s):=n(x,t_j+s),
	\qquad (x,s)\in\Omega\times (0,\tau),
  \eas
  and
  \be{30.88}
	a_j(s):=a(t_j+s),
	\qquad s\in (0,\tau),
  \ee
  for $j\in\N$. Then (\ref{30.8}) implies that
  \bas
	\int_0^\tau \io |n_j(x,s)-a_j(s)|^2 dxds
	&=& \int_{t_j}^{t_j+\tau} \io |n(x,t)-a(t)|^2 dxdt \\[3mm]
	&\to& 0
	\qquad \mbox{as } j\to\infty,
  \eas
  meaning that for
  \bas
	z_j(x,s):=n_j(x,s)-a_j(s),
	\qquad (x,s)\in \Omega\times (0,\tau), \qquad j\in\N,
  \eas
  we have
  \be{30.9}
	z_j\to 0
	\quad \mbox{in } L^2(\Omega\times (0,\tau))
	\qquad \mbox{as } j\to\infty.
  \ee
  Again by (\ref{30.22}), it follows from (\ref{30.66}) and (\ref{30.88}) that $(a_j)_{j\in\N}$ is bounded in
  $L^2((0,\tau))$, whence passing to a subsequence if necessary we may assume that for some nonnegative
  $a_\infty\in L^2((0,\tau))$ we have
  \be{30.10}
	a_j\wto a_\infty
	\quad \mbox{in } L^2((0,\tau))
	\qquad \mbox{as } j\to\infty.
  \ee
  Therefore, in view of (\ref{mass_n}) we see that
  \bas
	\int_0^\tau \io z_j(x,s)dxds
	&=& \int_0^\tau \io \Big(n_j(x,s)-a_j(s)\Big) dxds \\
	&=& \tau |\Omega| \onz - |\Omega| \cdot \int_0^\tau a_j(s)ds \\
	&\to& \tau |\Omega| \onz - |\Omega| \cdot \int_0^\tau a_\infty(s)ds 
	\qquad \mbox{as } j\to\infty,
 \eas
  which combined with (\ref{30.9}) allows us to determine the integral of the limit in (\ref{30.10}) according to
  \be{30.11}
	\int_0^\tau a_\infty(s)ds = \tau \onz.
  \ee
  On the other hand, rewriting (\ref{30.4}) in terms of $n_j$ and integrating in time we see that
  \bea{30.12}
	\frac{C_1 \tau}{4}
	&\le& \int_0^\tau \io n_j(x,s)\psi(x) dxds - \int_0^\tau \io \onz \psi(x) dxds \nn\\
	&=& \int_0^\tau \io n_j(x,s)\psi(x) dxds -\tau \onz \cdot \io \psi(x)dx
	\qquad \mbox{for all } j\in\N,
  \eea
  where as a consequence of (\ref{30.9}) and (\ref{30.10}),
  \bas
	\int_0^\tau \io n_j(x,s)\psi(x)dxds
	&=& \int_0^\tau \io z_j(x,s) \psi(x)dxds
	+ \int_0^\tau \io a_j(s) \psi(x)dxds \\
	&=& \int_0^\tau \io z_j(x,s) \psi(x)dxds
	+ \bigg(\int_0^\tau a_j(s) \bigg) \cdot \bigg( \io \psi(x)dx\bigg) \\
	&\to& 
	\bigg(\int_0^\tau a_\infty(s) \bigg) \cdot \bigg( \io \psi(x)dx\bigg) 
	\qquad \mbox{as } j\to\infty.
  \eas
  Taking $j\to\infty$ in (\ref{30.12}) we thus arrive at the conclusion that
  \bas
	\frac{C_1 \tau}{4}
	\le \bigg(\int_0^\tau a_\infty(s) \bigg) \cdot \bigg( \io \psi(x)dx\bigg) 
	-\tau \onz \cdot \io \psi(x) dx,
  \eas
  which in light of (\ref{30.11}) is absurd and hence proves that actually (\ref{30.1}) is valid.
\qed
{\bf Remark.} \quad
Let us mention here that in the case $m<2$ one may alternatively prove Lemma \ref{lem30} by 
invoking standard H\"older estimates for solutions of degenerate parabolic equations (\cite{porzio_vespri}): 
In fact, based on our previous estimates we can in this case obtain a bound for $n$, independent of $t>1$, 
in the space $C^{\theta,\frac{\theta}{2}}(\bar\Omega\times [t,t+1])$ for some $\theta>0$. 
This can be combined with the outcome of Lemma \ref{lem22} to see that actually $n(\cdot,t)\to\onz$
in $L^\infty(\Omega)$ as $t\to\infty$.\abs
We next make essential use of the fact that $f$ does not have positive zeroes to verify that $c$ decays
in the claimed sense.
\begin{lem}\label{lem29}
  Let $m>\frac{7}{6}$, and assume (\ref{f_pos}). Then the solution of (\ref{0}) constructed in Theorem \ref{theo16}
  satisfies
  \be{29.1}
	c(\cdot,t)\to 0 \quad \mbox{in } L^\infty(\Omega)
	\qquad \mbox{as } t\to\infty.
  \ee
\end{lem}
\proof
  If the claim was false, the for some $C_1>0$, some $(x_j)_{j\in\N}\subset\Omega$ and some 
  $(t_j)_{j\in\N} \subset (0,\infty)$ such that $t_j\to\infty$ as $j\to\infty$ we would have
  \bas
	c(x_j,t_j)\ge C_1
	\qquad \mbox{for all } j\in\N,
  \eas
  where passing to subsequences we may assume that there exists $x_0\in\bar\Omega$ such that
  $x_j\to\infty$ as $j\to\infty$. Since $c$ is uniformly continuous in 
  $\bigcup_{j\in\N} \Big(\bar\Omega \times [t_j,t_j+1]\Big)$
  by Lemma \ref{lem27}, this entails that on extracting a further subsequence if necessary, we can find $\delta>0$
  and $\tau\in (0,1)$ such that with $B:=B_\delta(x_0)\cap\Omega$ we have
  \bas
	c(x,t)\ge\frac{C_1}{2}
	\qquad \mbox{for all $x\in B, t\in (t_j,t_j+\tau)$ and } j\in\N,
  \eas
  so that since $f>0$ on $(0,\infty)$ by assumption (\ref{f_pos}), we see that
  \be{29.2}
	f(c(x,t)) \ge C_2
	\qquad \mbox{for all $x\in B, t\in (t_j,t_j+\tau)$ and } j\in\N
  \ee
  with some $C_2>0$.
  Now writing 
  \bas
	n_j(x,s):=n(x,t_j+s)
	\quad \mbox{and} \quad
	c_j(x,s):=c(x,t_j+s)
  \eas
  for $x\in\Omega, s\in (0,\tau)$ and $j\in\N$, from Lemma \ref{lem21} we obtain that
  \bea{29.3}
	\int_0^\tau \int_B n_j(x,s) f(c_j(x,s)) dxds
	&=& \int_{t_j}^{t_j+\tau} \int_B n(x,t) f(c(x,t)) dxdt \nn\\
	&\le& \int_{t_j}^\infty \int_B n(x,t) f(c(x,t)) dxdt \nn\\[2mm]
	&\to& 0
	\qquad \mbox{as } j\to\infty.
  \eea
  On the other hand, if we let $\psi(x):=\chi_B(x)$ for $x\in\Omega$, then from Lemma \ref{lem30} we obtain that
  \bas
	\bigg| \int_{t_j}^{t_j+\tau} \io n(x,t) \psi(x) dxdt - \nzb \tau |B| \bigg|
	&=& \Bigg| \int_{t_j}^{t_j+\tau} \bigg\{ \io n(x,t)\psi(x) dx - \io \onz \psi(x)dx \bigg\} \, dt \Bigg| \\
	&\le& \tau \cdot \sup_{t\in (t_j,t_j+\tau)} \bigg| \io n(x,t)\psi(x)dx - \io \nzb \psi(x)dx \bigg| \\[3mm]
	&\to& 0
	\qquad \mbox{as } j\to\infty,
  \eas
  and that hence
  \be{29.4}
	\int_0^\tau \int_B n_j(x,s)dxds \to \nzb \tau |B|
	\qquad \mbox{as } j\to\infty.
  \ee
  Therefore, (\ref{29.2}) warrants that the expression on the left of (\ref{29.3}) actually satisfies
  \bas
	\liminf_{j\to\infty}
	\int_0^\tau \int_B n_j(x,s) f(c_j(x,s))dxds
	&\ge& \liminf_{j\to\infty} \bigg\{ C_2 \int_0^\tau \int_B n_j(x,s)dxds \bigg\} \\[3mm]
	&=& C_2 \cdot \nzb \tau |B|,
  \eas
  which contradicts (\ref{29.3}) and hence proves the lemma.
\qed
Finally, decay of $u$ will be a consequence of the stabilization property of $n$ in Lemma \ref{lem30}.
Since the latter has been asserted in the weak-$\star$ sense in $L^\infty(\Omega)$ only, an argument
based on the use of appropriate linear functionals involving $u$ seems adequate to derive this here.
\begin{lem}\label{lem31}
  Let $m>\frac{7}{6}$. Then with $(n,c,u)$ as given by Theorem \ref{theo16}, we have
  \be{31.1}
	u(\cdot,t) \to 0
	\quad \mbox{in } L^\infty(\Omega)
	\qquad \mbox{as } t\to\infty.
  \ee
\end{lem}
\proof
  Since the Hilbert space realization $A_2$ of the Stokes operator is positive and self-adjoint in $L^2_\sigma(\Omega)$ 
  with compact inverse, there exists a complete
  orthonormal basis $(\psi_k)_{k\in\N}$ of eigenfunctions $\psi_k$ of $A_2$ corresponding to positive eigenvalues
  $\lambda_k$, $k\in\N$. 
  By density of $\bigcup_{N\in\N} {\rm span} \{\psi_k \ | \ k\le N\}$ in $L^2_\sigma(\Omega)$, in view of the 
  uniform H\"older continuity of $u$ in $\Omega\times (0,\infty)$, as asserted by Lemma \ref{lem28} and Lemma
  \ref{lem26}, to prove (\ref{31.1}) it is sufficient to show that for each $k\in\N$ we have
  \be{31.2}
	\io u(x,t) \cdot \psi_k(x)dx \to 0
	\qquad \mbox{as } k\to\infty.
  \ee
  For this purpose, we fix any such $k$ and let
  \bas
	y_\eps(t):=\io \ueps(x,t)\cdot\psi_k(x)dx, \qquad t\ge 0,
  \eas
  for $\eps\in (0,1)$. Then from (\ref{0eps}) and the eigenfunction property of $\psi_k$ we obtain that since
  $\nabla\cdot\psi_k\equiv 0$ and $\io \neps(\cdot,t)\equiv \onz$ by (\ref{mass}), 
  \bas
	y_\eps'(t)
	&=& - \io A \ueps \cdot \psi_k + \io \neps \nabla \phi \cdot\psi_k \\
	&=& - \io \ueps \cdot A\psi_k + \io (\neps-\onz) \nabla\phi \cdot\psi_k + \onz \io \nabla \phi \cdot\psi_k \\
	&=& - \lambda_k \io \ueps \cdot \psi_k + \io (\neps-\onz) \nabla \phi\cdot\psi_k \\[1mm]
	&=& - \lambda_k y_\eps(t) + g_\eps(t)
	\qquad \mbox{for all } t>0
  \eas
  with
  \bas
	g_\eps(t):=\io \Big(\neps(x,t)-\onz\Big) \nabla\phi \cdot \psi_k(x)dx,
	\qquad t\ge 0.
  \eas
  Upon integration, this shows that for any choice of $t_0\ge 0$,
  \bas
	y_\eps(t)=y_\eps(t_0) e^{-\lambda_k(t-t_0)} + \int_{t_0}^t e^{-\lambda_k(t-s)} g_\eps(s) ds
	\qquad \mbox{for all } t>t_0.
  \eas
  Since as $\eps=\eps_j\searrow 0$ we have 
  $\ueps \to u$ in $C^0_{loc}(\bar\Omega\times [0,\infty))$ and $\neps\wsto n$ in $L^\infty(\Omega\times (0,\infty))$
  by Lemma \ref{lem26}, we may take $\eps=\eps_j \searrow 0$ here to infer that with
  \bas
	y(t):=\io u(x,t) \cdot \psi_k(x)dx
	\quad \mbox{and} \quad
	g(t):=\io \Big(n(x,t)-\onz\Big) \nabla\phi\cdot\psi_k(x)dx,
	\qquad t\ge 0,
  \eas
  we have
  \be{31.3}
	y(t)=y(t_0) e^{-\lambda_k(t-t_0)} + \int_{t_0}^t e^{-\lambda_k(t-s)} g(s)ds
	\qquad \mbox{for all $t_0\ge 0$ and any } t>t_0,
  \ee
  where thanks to Lemma \ref{lem30} we know that
  \bas
	g(t)\to 0
	\qquad \mbox{as } t\to\infty.
  \eas
  Accordingly, if in order to prove (\ref{31.2}) we let $\delta>0$ be given, then we can pick $t_0>0$ large enough
  fulfilling
  \be{31.4}
	|g(t)| < \frac{\lambda_k \delta}{2}
	\qquad \mbox{for all } t>t_0.
  \ee
  As $C_1:=\|y\|_{L^\infty((0,\infty))}$ is finite due to the boundedness of $u$ in $\Omega\times (0,\infty)$ guaranteed
  by Lemma \ref{lem15} and Lemma \ref{lem26}, from (\ref{31.3}) we thus infer that
  \bas
	|y(t)|
	&\le& C_1 e^{-\lambda_k(t-t_0)}
	+ \frac{\lambda_k \delta}{2} \int_{t_0}^t e^{-\lambda_k(t-s)} ds \\
	&=& C_1 e^{-\lambda_k(t-t_0)}
	+ \frac{\lambda_k \delta}{2} \cdot \frac{1-e^{-\lambda_k(t-t_0)}}{\lambda_k} \\
	&<& C_1 e^{-\lambda_k(t-t_0)}
	+ \frac{\delta}{2}
	\qquad \mbox{for all } t>t_0.
  \eas
  This implies that with $t_1:=\max\{t_0,t_0+\frac{1}{\lambda_k} \ln \frac{2C_1}{\delta}\}$ we have
  \bas
	|y(t)| < \delta
	\qquad \mbox{for all } t>t_1,
  \eas
  which establishes (\ref{31.2}) and thereby completes the proof.
\qed
\proofc of Theorem \ref{theo166}.\quad
  The claimed convergence properties are precisely asserted by Lemma \ref{lem30}, Lemma \ref{lem29} and Lemma \ref{lem31}.
\qed

\end{document}